\documentclass[12pt]{article}
\usepackage{amscd, amssymb, latexsym, amsmath, amsthm}

\newtheorem{thm}{Theorem}
\newtheorem{cor}{Corollary}
\newtheorem{pro}{Proposition}
\newtheorem{lem}{Lemma}
\newtheorem{dfn}{Definition}
\newtheorem{exa}{Example}
\newtheorem{rem}{Remark}
\newtheorem*{ack}{Acknowledgements}

\numberwithin{thm}{section} \numberwithin{cor}{section}
\numberwithin{pro}{section} \numberwithin{lem}{section}
\numberwithin{dfn}{section}
\numberwithin{rem}{section} \numberwithin{equation}{section}

\newcommand{\R}{\mathbb R}

\newcommand{\heat}{\left(\partial_t-\frac{1}{n+2} \Delta\right)}

\newcommand{\pt}{\partial_t}
\newcommand{\sfrac}[2]{{\textstyle \frac{#1}{#2}}}

\newcommand{\sph}{\mathbb S}

\begin{document}
\title{Ancient Solutions of the Affine Normal Flow}
\author{John Loftin\footnote{The first author is
partially supported by NSF Grant DMS0405873 and by the IMS at the
Chinese University of Hong Kong, where, in November, 2004, he
developed and presented some of the material in Sections
\ref{aff-nor-sec} and \ref{aff-str-sec} to students there.} \,and
Mao-Pei Tsui} \maketitle

\section{Introduction}

Consider a smooth, strictly convex hypersurface $\mathcal L$
locally parametrized by $F(x)\in\R^{n+1}$. The affine normal is a
vector field $\xi=\xi_{\mathcal L}$ transverse to $\mathcal L$ and
invariant under volume-preserving affine transformations of
$\R^{n+1}$. The affine normal flow evolves such a hypersurface in
time $t$ by
$$ \partial_t F(x,t) = \xi(x,t), \qquad F(x,0)=F(x).$$

In \cite{chow85}, Ben Chow proved that every smooth, strictly
convex hypersurface in $\R^{n+1}$ converges in finite time under
the affine normal flow to a point. In \cite{andrews96}, Ben
Andrews proved that the rescaled limit of the contracting
hypersurface around the final point converges to an ellipsoid.
Later, Andrews \cite{andrews00} also studied the case in which the
initial hypersurface is compact and convex with no regularity
assumed.  In this case, the affine normal flow, unlike the Gauss
curvature flow, is instantaneously smoothing. In other words, such
an initial hypersurface under the affine normal flow will evolve
to be smooth and strictly convex at any positive time before the
extinction time.

In the present work, we develop the affine normal flow for any
noncompact convex hypersurface $\mathcal L$ in $\R^{n+1}$ whose
convex hull $\hat{\mathcal L}$ contains no lines (if $\hat{\mathcal
L}$ contains a line, the affine normal flow does not move it at
all). As in \cite{andrews00} we define the flow by treating the
$\mathcal L$ as a limit of a nested sequence of  smooth, compact,
strictly convex hypersurfaces $\mathcal L^i$.  Our main new result
is to classify ancient solutions---solutions defined for time
$(-\infty,T)$---for the affine normal flow.
\begin{thm} \label{classify-ancient}
Any ancient solution to the affine normal flow must be be either an
elliptic paraboloid (which is a translating soliton) or an ellipsoid
(which is a shrinking soliton).
\end{thm}

The proof of Theorem \ref{classify-ancient} relies on a decay
estimate of Andrews for the cubic form $C^i_{jk}$ of a compact
hypersurface under the affine normal flow \cite{andrews96}. In
particular, the norm squared $|C|^2$ of the cubic form with respect
to the affine metric decays like $1/t$ from the initial time.  For
an ancient solution then, we may shift the initial time as far back
as we like, and thus the cubic form $C^i_{jk}$ is identically zero.
Then a classical theorem of Berwald shows that the hypersurface must
be a hyperquadric, and the paraboloid and ellipsoid are the only
hyperquadrics which form ancient solutions to the affine normal flow
(the hyperboloid, an expanding soliton, is not part of an ancient
solution).

In order to apply this estimate in our case, we need local
regularity estimates to ensure that for all positive time $t$, the
evolving hypersurfaces $\mathcal L^i(t)$ converge locally in the
$C^\infty$ topology to $\mathcal L(t)$.  Thus Andrews's pointwise
bound on the cubic form survives in the limit.  We work in terms
of the support function.  The $C^2$ estimates are provided by a
speed bound of Andrews \cite{andrews00} and a Pogorelov-type
Hessian bound similar similar to one in Guti\'errez-Huang
\cite{gutierrez-huang98}. These estimates provide uniform local
parabolicity, and then Krylov's theory and standard bootstrapping
provide local estimates to any order.

Another key ingredient is the use of barriers.  Here the invariance
of the affine normal flow under volume-preserving affine
transformations is important.  The main barriers we use are
ellipsoids and a particular expanding soliton (a hyperbolic affine
sphere) due to Calabi \cite{calabi72}.  In particular,
Guti\'errez-Huang's estimate can only be applied to solutions of
PDEs which move in time by some definite amount.  Calabi's example
is a crucial element in constructing a barrier to guarantee the
solution does not remain constant in time.

Solitons of the affine normal flow have been very well studied
\cite{calabi72,cheng-yau86}. They are precisely the affine spheres.
The shrinking solitons of the affine normal flow are the elliptic
affine spheres, and Cheng-Yau proved that any properly embedded
elliptic affine sphere must be an ellipsoid \cite{cheng-yau86}.
Translating solitons are parabolic affine spheres, and again
Cheng-Yau showed that any properly embedded parabolic affine sphere
must be an elliptic paraboloid \cite{cheng-yau86}.

Expanding solitons are hyperbolic affine spheres, which behave quite
differently. Cheng-Yau proved that every convex cone in $\R^{n+1}$
which contains no lines admits a unique (up to scaling) hyperbolic
affine sphere which is asymptotic to the boundary of the cone
\cite{cheng-yau77,cheng-yau86}. (For example, the hyperboloid is the
hyperbolic affine sphere asymptotic to the standard round cone.) The
converse is also true: every properly embedded hyperbolic affine
sphere in $\R^{n+1}$ is asymptotic to the boundary of a convex cone
containing no lines \cite{cheng-yau86}.  Our definition of the
affine normal flow immediately provides an expanding soliton which
is a weak (viscosity) solution, and our local regularity estimates
show that this solution is smooth.

We should note that Cheng-Yau \cite{cheng-yau86} proved results for
hyperbolic affine sphere based on the affine metric.  In particular,
a hyperbolic affine sphere has complete affine metric if and only if
it is properly embedded in $\R^{n+1}$ if and only if it is
asymptotic to the boundary of a convex cone in $\R^{n+1}$ containing
no lines. Our methods do not yet yield any insight into the affine
metric of evolving hypersurfaces.  If the initial hypersurface of
the affine normal flow is the boundary of a convex cone containing
no lines, then at any positive time, the solution is the
homothetically expanding hyperbolic affine sphere asymptotic to the
cone.  Cheng-Yau's result implies the affine metric in this case is
complete at any positive time $t$.  It will be interesting to
determine whether, under the affine normal flow, the affine metric
is complete at any positive time for any noncompact properly
embedded initial hypersurface.  Presumably a parabolic version of
the affine geometric gradient estimate of Cheng-Yau is needed, as
suggested by Yau \cite{yau05}.

When restricted to an affine hyperplane, the support function of a
hypersurface evolving under the affine normal flow satisfies
\begin{equation}
  \label{our-eq}
\partial_t s = -(\det\partial^2_{ij}s)^{-\frac1{n+2}}.
\end{equation}
Guti\'errez and Huang \cite{gutierrez-huang98} have studied a
similar parabolic Monge-Amp\`ere equation
$$ \partial_t s = - (\det \partial^2_{ij}s)^{-1}. $$
They prove that any ancient entire solution to this equation which
\emph{a priori} satisfies bounds on the ellipticity must be an
evolving quadratic polynomial.  Our Theorem \ref{classify-ancient}
reduces to an similar result for (\ref{our-eq}):  The ellipsoid and
paraboloid solitons provide ancient solutions to (\ref{our-eq})
which can be represented, up to possible affine coordinate changes,
by
$$ s=\left(-\frac{2n+2}{n+2}\,t\right)^{\frac{n+2}{2n+2}} \sqrt{1+|y|^2},
\qquad s = \frac{|y|^2}2 - t$$ respectively. Our result doesn't
require any \emph{a priori} bounds on the ellipticity. We do not
require our solutions to be entire, but they do solve a Dirichlet
boundary condition. See Section \ref{dirichlet-section} below.

We also mention a related theorem due to J\"orgens \cite{jorgens54}
for $n=2$, Calabi \cite{calabi58} for $n\le 5$, and independently to
Pogorelov \cite{pogorelov72} and Cheng-Yau \cite{cheng-yau86} for
all dimensions:
\begin{thm}
Any entire convex solution to $$\det \partial^2_{ij} u = c>0$$ is an
quadratic polynomial.
\end{thm}
The graph of each such $u$ is a parabolic affine sphere, and
Cheng-Yau's classification provides the result. Our techniques do
not yet yield an independent proof of this classical theorem: We do
not yet know if the affine normal flow is unique for a given initial
convex noncompact hypersurface. Even though any parabolic affine
sphere may naturally be thought of as a translating soliton under
the affine normal flow, the flow we define, with the parabolic
affine sphere as initial condition, may not \emph{a priori} be the
same flow as the soliton solution, and thus may not come from an
ancient solution in our sense.

It is also interesting to compare our noncompact affine normal flow
with other geometric flows on noncompact hypersurfaces.  In
particular, Ecker-Huisken and Ecker have studied mean-curvature flow
of entire graphs in Euclidean space \cite{ecker-huisken89}
\cite{ecker-huisken91} and of spacelike hypersurfaces in Lorentzian
manifolds  \cite{ecker93} \cite{ecker97} \cite{ecker03}.  In
\cite{ecker-huisken91}, Ecker-Huisken prove that under any entire
graph of a locally Lipschitz function moves under the mean curvature
flow in Euclidean space to be smooth at any positive time, and the
solution exists for all time.  Ecker proves long-time existence for
any initial spacelike hypersurface in Minkowski space under the mean
curvature flow \cite{ecker97} and proves instantaneous smoothing for
some weakly spacelike hypersurfaces in \cite{ecker03}.

In the present work, we prove instantaneous smoothing and
long-time existence for the affine normal flow on noncompact
hypersurfaces for any initial convex noncompact properly embedded
hypersurface $\mathcal L\subset \R^{n+1}$  which contains no
lines.  In this case, the evolving hypersurface $\mathcal L(t)$
under the affine normal flow exists for all time $t>0$ (Theorem
\ref{long-time-ex}) and is smooth for all $t>0$ (Theorem
\ref{smooth-hypersurface}).  Moreover, the following maximum
principle at infinity is satisfied: If $\mathcal L^1$ and
$\mathcal L^2$ are convex properly embedded hypersurfaces whose
convex hulls  satisfy $\widehat {\mathcal L^1} \subset
\widehat{\mathcal L^2}$, then for all $t>0$, the convex hulls
satisfy $\widehat{\mathcal L^1(t)} \subset \widehat {\mathcal
L^2(t)}$. This sort of maximum principle at infinity does not hold
for all evolution equations of noncompact hypersurfaces.  In
particular, there is an example due to Ecker \cite{ecker97}, of
two soliton solutions to the mean curvature flow in Minkowski
space, for which this fails.

The affine normal flow is equivalent (up to a diffeomorphism) to
the hypersurface flow by $K^{\frac1{n+2}}\nu$, where $K$ is the
Gauss curvature and $\nu$ is the inward unit normal.  The
techniques we use (the definitions and ellipticity estimates)
should apply to flows of noncompact convex hypersurfaces by other
power of the Gauss curvature.  Andrews \cite{andrews00} addresses
many aspects of the compact case of flow by powers of Gauss
curvature.  In particular, he verifies that for $\alpha\le 1/n$,
any convex compact hypersurface in $\R^{n+1}$ evolves under the
flow by $K^\alpha \nu$ to be smooth and strictly convex at any
positive time $t$.  In essence, we verify this in the noncompact
case for $\alpha=1/(n+2)$ (see Theorem \ref{smooth-hypersurface}
below).  We expect the same result to be true in the noncompact
case for all $\alpha\le 1/n$. We should note that for
$\alpha>1/n$, flat sides of any initial hypersurface remain
non-strictly convex for some positive time. We note that in the
case of the Gauss curvature flow in $\R^3$ ($\alpha=1$),
Daskalopoulos-Hamilton \cite{daskalopoulos-hamilton99} study how
the boundary of such a flat side evolves over time.

Our treatment of the affine normal flow is largely self-contained.
In Sections \ref{aff-nor-sec} and \ref{aff-str-sec}, we recall the
definition of the affine normal and the basic affine structure
equations.  We develop the computations necessary by using
notation similar to that of e.g.\ Zhu \cite{zhu02}: let
$F\!:U\to\R^{n+1}$ represent a local embedding of a hypersurface
for $U\subset\R^n$ a domain.  Then we derive the structure
equations based on derivatives of $F$. Using this notation, we
develop the affine normal flow of the basic quantities associated
with the hypersurface in Sections \ref{flow-euc-sec},
\ref{flow-C-sec} and \ref{flow-A-sec}.  The main estimate we need
on the cubic form is found in Section \ref{flow-C-sec}. These
evolution equations are all due to Andrews \cite{andrews96}, and
we include derivations of them for the reader's convenience.  In
Section \ref{support-function-section}, we introduce the support
function and some basic results we will need. We define our affine
normal flow on a noncompact convex hypersurface $\mathcal L$ in
Section \ref{flow-def-sec}, basically as a limit of compact convex
hypersurfaces approaching $\mathcal L$ from the inside, and we
verify that the soliton solutions behave properly under our
definition in Section \ref{sol-ex-sec}.

In Section \ref{andrews-est-sec}, we turn to the estimates that
are the technical heart of the paper.  We prove an estimate of
Andrews on the speed of the support function evolving under affine
normal flow \cite{andrews00}.  In particular, we verify that this
estimate survives in the limit to our noncompact hypersurface.  In
Section \ref{gh-est-sec}, we prove a version of a Pogorelov-type
estimate due to Guti\'errez-Huang \cite{gutierrez-huang98}, which
bounds the Hessian of the evolving support function, and in
Section \ref{upper-interior-barrier}, we construct barriers to
ensure that Guti\'errez-Huang's estimate applies. Krylov's
estimates then ensure the support function is smooth for all time
$t>0$. In Section \ref{hyp-reg-sec}, we verify that the evolving
hypersurface is smooth as well, and relate the noncompact affine
normal flow to a Dirichlet problem for the support function in
Section \ref{dirichlet-section}.  The main results are proved in
Section \ref{final-sec}.

Our treatment of noncompact hypersurfaces as limits of compact
hypersurfaces is a bit different from the usual analysis on
noncompact manifolds. Typically noncompact manifolds are exhausted
by compact domains with boundary (e.g.\ geodesic balls on complete
Riemannian manifolds or sublevel sets of a proper height function
on a hypersurface considered as a Euclidean graph), and then a
version of the maximum principle is shown to hold in the limit of
the exhaustion.  Our limiting process is extrinsic, on the other
hand: We apply the maximum principle to $|C|^2$ to derive
Andrews's pointwise bound on compact hypersurfaces without
boundary, which in turn survives in the limiting noncompact
hypersurface. It is still desirable to implement an approach by
intrinsically exhausting the hypersurface, to be able to use the
maximum principle more directly on the evolving noncompact
hypersurface. Perhaps the description in Section
\ref{dirichlet-section} of the affine normal flow in terms of a
Dirichlet problem for the support function will be of some use.

\begin{ack}
We would like to thank S.T.\ Yau for introducing us to the beautiful
theory of affine differential geometry, Richard Hamilton for many
inspiring lectures on geometric evolution equations, and D.H.\ Phong
for his constant encouragement.
\end{ack}

Notation: Subscripts after a comma are used to denote covariant
derivatives with respect to the affine metric.  So the second
covariant derivative of $H$ is $H_{,ij}$, for example. Of course the
first covariant derivative of a function is just ordinary
differentiation, which commutes with the time derivative
$\partial_t$. $\partial_i$ will denote an ordinary space derivative.
We use Einstein's summation convention that any paired indices, one
up and one down, are to be summed from 1 to $n$. Unless otherwise
noted, we raise and lower indices using the affine metric $g_{ij}$.

\section{The affine normal} \label{aff-nor-sec}

Here we define the affine normal to a hypersurface in a similar way
to Nomizu-Sasaki \cite{nomizu-sasaki}, but using notation adapted to
our purposes.

Let $F=F(x^1,\dots,x^n)$ be a local embedding of a smooth,
strictly convex hypersurface in $\R^{n+1}$. Let $F\!: \Omega\to
\R^{n+1}$, where $\Omega$ is a domain in $\R^n$. Let $\tilde\xi$
be a smooth transverse vector field to $F$.  Now we may
differentiate to determine
\begin{eqnarray} \label{2-d-f}
\partial^2_{ij} F &=& \tilde g_{ij}\tilde \xi + \tilde \Gamma^k_{ij} \partial_k F,\\
\label{space-deriv-xi}
\partial_i \tilde \xi &=& \tilde \tau_i \tilde \xi - \tilde A^j_i \partial_j F.
\end{eqnarray}
It is straightforward to check that $\tilde g_{ij}$ is a symmetric
tensor, $\tilde \Gamma^k_{ij}$ is a torsion free connection,
$\tilde\tau_i$ is a one-form, and $\tilde A^j_i$ is an
endomorphism of the tangent bundle.  With respect to $\tilde \xi$,
$\tilde g_{ij}$ is called the \emph{second fundamental form} and
$\tilde A^j_i$ is the \emph{shape operator}.

\begin{pro} \label{aff-normal-def}
There is a unique transverse vector field $\xi$,
called the \emph{affine normal}, which satisfies
\begin{enumerate}
\item \label{inward} $\xi$ points inward. In other words, $ \xi$
and the hypersurface $F(\Omega)$ are on the same side of the
tangent plane.
 \item \label{equiaffine} $\tau_i = 0$.
 \item \label{vol-form}
 $\det  g_{ij} = \det (\partial_1 F,\dots,
\partial_n F, \xi)^2.$  The determinant on the left is that of an
 $n\times n$ matrix, while the determinant on the right is that on
 $R^{n+1}$.
\end{enumerate}
\end{pro}

Note we have dropped the tildes in quantities defined by the
affine normal (the connection term is an exception: see the next
section). Condition \ref{inward} implies that the second
fundamental form $g_{ij}$ is positive definite, and thus we say
$g_{ij}$ is the \emph{affine metric}. Condition \ref{equiaffine}
is called that $\xi$ is \emph{equiaffine}. Condition
\ref{vol-form} is that the volume form on the hypersurface induced
by $\xi$ and the volume form on $\R^{n+1}$ is the same as the
volume form induced by the affine metric.

The following proof of Proposition \ref{aff-normal-def} will be
instructive in computing the affine normal later on.

\begin{proof}
Given an arbitrary inward-pointing transverse vector filed $\tilde
\xi$, any other may be written as $\xi = \phi \tilde \xi + Z^i
\partial_i F$, where $\phi$ is a positive scalar function and
$Z^i\partial_i F$ is a tangent vector field.

Condition \ref{vol-form} determines $\phi$ in terms of $\tilde
\xi$:  Plug  $\xi = \phi\tilde\xi + Z^i \partial_i F$ into
(\ref{2-d-f}), and the terms in the span of $\tilde \xi$ give
\begin{equation} \label{affine-met-def}
g_{ij} = \phi^{-1}\tilde g_{ij}.
\end{equation}
Now Condition \ref{vol-form} shows that
$$ \phi^{-n} \det \tilde g_{ij} = \det g_{ij} = \det(\partial_1 F,
\dots , \partial_n F, \xi)^2 = \phi^2 \det(\partial_1 F, \dots,
\partial_n F,\tilde \xi)^2,$$
and so
\begin{equation} \label{phi-def}
\phi = \left( \frac{ \det \tilde g_{ij}} {\det (\partial_1 F,
\dots, \partial_n F, \tilde \xi)^2} \right)^\frac1{n+2}.
\end{equation}

Finally, we use the equiaffine condition to determine $Z^i$: Plug
in for $\xi$, set $\tau_i=0$, and consider the terms in the span
of $\tilde\xi$ to find
\begin{eqnarray}
\nonumber
-A_i^j\partial_j F &=&\partial_i (\phi\tilde \xi + Z^j
\partial_j F) \\
\nonumber
&=& \partial_i \phi\,\tilde \xi + \phi\,
\partial_i\tilde\xi +
\partial_i Z^j \,\partial_j F + Z^j \,\partial^2_{ij} F \\
\nonumber
&=& \partial_i \phi\,\tilde \xi + \phi (\tilde\tau_i
\tilde\xi - \tilde A^j_i \partial_j F) +
\partial_i Z^j \,\partial_j F + Z^j (\tilde g_{ij} \xi + \tilde
\Gamma^k_{ij} \partial_k F),\\
\nonumber
0&=& \partial_i\phi + \phi\,\tilde\tau_i + Z^j\tilde g_{ij}, \\
\label{Z-def} Z^j &=& - \tilde g^{ij} (\partial_i \phi +
\phi\,\tilde\tau_i),
\end{eqnarray}
where $\tilde g^{ij}$ is the inverse matrix of $\tilde g_{ij}$.
\end{proof}

\begin{cor}
The affine normal is invariant under volume-preserving affine
automorphisms of $\R^{n+1}$.  In other words, if $\Phi$ is such an
affine map, and $\xi$ is the affine normal filed to a hypersurface
$F(\Omega)$, then $\Phi_*\xi$ is the affine normal to $(\Phi\circ
F)(\Omega)$.
\end{cor}

\begin{proof}
The defining conditions in the proposition are invariant under
affine volume-preserving maps on $\R^{n+1}$.
\end{proof}

\section{Affine structure equations} \label{aff-str-sec}

Consider a smooth, strictly convex hypersurface in $\R^{n+1}$
given by the image of an embedding $F=F(x^1,\dots,x^n)$.  The
affine normal is an inward-pointing transverse vector field to the
hypersurface, and we have the following structure equations:
\begin{eqnarray}
\partial^2_{ij} F &=& g_{ij} \xi + (\Gamma^k_{ij} + C^k_{ij})
F_{,k} \label{s1} \\
\label{s2} \xi_{,i} &=& -A^k_i F_{,k}
\end{eqnarray}
Here $g_{ij}$ is the \emph{affine metric}, which is positive
definite. $\Gamma^k_{ij}$ are the Christoffel symbols of the
metric. Since $\Gamma^k_{ij} + C^k_{ij}$ is a connection, then
$C^k_{ij}$ is a tensor called the \emph{cubic form}. $A^k_i$ is
the \emph{affine curvature}, or \emph{affine shape operator}.
Equation (\ref{s1}) shows immediately that $$ C^k_{ij} =
C^k_{ji}.$$

Now consider the second covariant derivatives with respect to the
affine metric
\begin{eqnarray} \nonumber
F_{,ij} &=& \partial^2_{ij} F - \Gamma^k_{ij} F_{,k} \\
 \label{F2cov} &=& g_{ij}\xi + C^k_{ij} F_{,k} \\
 \nonumber \xi_{,ij} &=& -A^k_{i,j} F_{,k} - A^k_i F_{,kj} \\
 \nonumber &=&-A^k_{i,j} F_{,k} - A_{ij}\xi - A^k_i C^\ell_{kj} F_{,\ell}
\end{eqnarray}
Since $\xi_{,ij} = \xi_{,ji}$, we have $$A_{ij}=A_{ji}$$ and the
following Codazzi equation for the affine curvature:
 \begin{equation} \label{A-cod}
 A_{j,i}^k -
A_{i,j}^k = A^\ell_i C^k_{\ell j} - A^\ell_j C^k_{\ell i},
\end{equation}
$$A_{jk,i}=A_{ji,k} + A^l_i C_{ljk} - A^l_j C_{lik}.$$
Finally, consider the third covariant derivative of $F$
\begin{eqnarray*}
F_{,ijk} &=& g_{ij} \xi_{,k} + C^\ell_{ij,k} F_\ell + C^\ell_{ij}
F_{,\ell k} \\
&=& - g_{ij} A_k^\ell F_{,\ell} + C^\ell_{ij,k} F_{,\ell} +
C_{ijk} \xi + C^m_{ij} C_{m k}^\ell F_{,\ell}
\end{eqnarray*}

Recall the conventions for commuting covariant derivatives of
tensors by using the Riemannian curvature $R^\ell_{ijk}$:
$$v^h_{,ji}-v^h_{,ij}=R^h_{ijk} v^k, \quad \mbox{and}\quad
w_{k,ji}-w_{k,ij}=-R^h_{ijk} w_h.$$ Therefore,
\begin{eqnarray*}
-R_{jki}^\ell F_{,\ell} &=& F_{,ikj} - F_{,ijk} \\
&=& - g_{ik} A_j^\ell F_{,\ell} + C^\ell_{ik,j} F_{,\ell} +
C_{ikj} \xi + C^m_{ik} C_{m j}^\ell F_{,\ell} \\
&& {} + g_{ij} A_k^\ell F_{,\ell} - C^\ell_{ij,k} F_{,\ell} -
C_{ijk} \xi - C^m_{ij} C_{m k}^\ell F_{,\ell}
\end{eqnarray*}
From the part of this equation in the span of $\xi$, we see
$$ C_{ikj} = C_{ijk},$$
and so the cubic form is totally symmetric in all three indices.
Lower the index $R_{jk\ell i} = R^m_{jki}g_{m\ell}$ and compute
$2R_{jk\ell i} = R_{jk\ell i}-R_{jk i\ell }$ to find
\begin{equation} \label{curv-eq}
R_{jk\ell i} =\sfrac12 g_{ik} A_{j\ell} - \sfrac12 g_{ij} A_{k\ell}
- \sfrac12 g_{\ell k} A_{ji} +\sfrac12 g_{\ell j} A_{ki} - C_{ik}^m
C_{mj\ell} + C^m_{ij} C_{mk\ell},
\end{equation}
$$R^\ell_{jki}= \sfrac12( g_{ik} A^{\ell}_j -g_{ij} A^{\ell}_k -\delta^{\ell}_k
A_{ji}+\delta^{\ell}_j A_{ki}) -
C^m_{ik}C^{\ell}_{mj}+C^m_{ij}C^{\ell}_{mk},$$ and the Ricci
curvature of the affine metric
$$R_{ki}=g^{j \ell}R_{jk\ell i}=\sfrac12 g_{ik} H+\sfrac{n-2}2 A_{ki}  + C^{m\ell}_{i}
C_{mk\ell}.$$ Note here that $H=A_i^i$ is the \emph{affine mean
curvature}.

On the other hand we may compute $0 = R_{jki\ell} + R_{jk\ell i}$ to
find the following Codazzi equation for the cubic form:
 \begin{equation} \label{C-cod}
C_{ij\ell,k} - C_{ik\ell,j} = \sfrac12 g_{ij}A_{k\ell} -\sfrac12
g_{ik} A_{j\ell} + \sfrac12 g_{\ell j} A_{ki} - \sfrac12 g_{\ell k}
A_{ji}.
\end{equation}

Thus far, we have only used equations (\ref{s1}) and (\ref{s2}) to
derive the structure equations.  The only constraint is that the
transversal vector field $\xi$ be equiaffine.  The position vector
and the Euclidean normal are also equiaffine.  Another important
property of the affine normal is the following \emph{apolarity}
condition
\begin{equation}\label{apolarity}
  C^i_{ij}=0,
\end{equation}
which follows from taking the first covariant derivative of
Condition \ref{vol-form} in Proposition \ref{aff-normal-def}:
\begin{eqnarray*}
0&=& \partial_j \det (F_{,1},\dots,F_{,n}, \xi ) \\
&=& \det (F_{,1j},\dots,F_{,n},\xi) + \cdots + \det (F_{,1},
\dots,
F_{,nj},\xi) + \det (F_{,1},\dots, F_{,n}, \xi_{,j}) \\
&=& C^1_{1j}\det (F_{,1},\dots,F_{,n}, \xi ) + \cdots + C^n_{nj}
\det (F_{,1},\dots,F_{,n}, \xi ) + 0 \\
&=& \left(C^i_{ij} \right) \det (F_{,1},\dots,F_{,n}, \xi ).
\end{eqnarray*}

The apolarity condition and (\ref{F2cov}) imply the following
formula for the affine normal in terms of the metric:
$$ \xi = \frac{\Delta F}n.$$

\section{Evolution of $\bar{g}_{ij}$, $g_{ij}$ and $K$}
\label{flow-euc-sec}
 Let $M^n$ be an $n$-dimensional smooth manifold
and let $F(\cdot, t) : M^n  \mapsto R^{n+1}$ be a one-parameter
family of smooth hypersurface immersions in $R^{n+1}$. We say that
it is a solution of the affine normal  flow if
\begin{eqnarray}
\label{anf}\partial_t F=\frac{\partial F(x, t)}{\partial t}= \xi \;
,\quad  x \in M^n \; , \quad t
> 0\end{eqnarray}
 where $\xi$  is affine normal  flow on $F(\cdot , t)$.

In a local coordinate system $\{x_i\},$ $1 \le  i \le n$. The
Euclidean inner product $\langle\cdot,\cdot\rangle$ on $\R^{n+1}$
induces the metric $\bar{g}_{ij}$ and the Euclidean second
fundamental form $h_{ij}$ on $F(\cdot,t)$. These can be computed as
follows
$$\bar{g}_{ij}= \langle \partial_i F, \partial_j
F \rangle$$ and
$$h_{ij}= \langle \partial^2_{ij} F, \nu
\rangle,$$ where $\nu$ is the unit inward normal on $F(\cdot, t)$.
The Gaussian curvature is
$$K=\frac{\det h_{ij}}{\det\bar{g}_{ij}}.$$
 By (\ref{affine-met-def}) and (\ref{phi-def}), the affine metric
is
$${g}_{ij}= \frac{h_{ij}}{\phi}, \quad \text{where} \quad \phi =
K^{\frac{1}{n+2}}.$$ (Note that $\det\bar g_{ij} = \det
(\partial_1F,\dots,\partial_n F,\nu)^2.$) Proposition
\ref{aff-normal-def} shows that the affine normal is
\begin{equation}
\label{aff-euc-normal} \xi= - h^{ki}\,\partial_i \phi \,\partial_k F
+\phi \nu= - g^{ki}\,\partial_i( \ln \phi) \,\partial_k F +\phi \nu.
 \end{equation}
(Note that $\nu$ is equiaffine.) Also recall the affine curvature
$\{A^k_j\}$ is defined by
 \begin{equation} \label{afc} \partial_j \xi
= -A^k_j \,\partial_k F.
\end{equation}

As we'll see below in Section \ref{support-function-section}, the
support function of a smooth convex hypersurface is defined by
$$s = -\langle F , \nu \rangle.$$

\begin{pro} \label{euclidean-evolve} Under the
affine normal  flow,
\begin{eqnarray*}
 \partial_t F_{,i}
&=&  -A^k_i F_{,k},\\
 \partial_t \nu &=& 0 ,\\
  \partial_j \nu &=& - h_{jl}
\bar{g}^{lm} F_{,m} ,\\
\partial_t\bar{g}_{ij} &=& - (A^k_i \bar{g}_{kj}
+A^k_j \bar{g}_{ki} ),\\
 \partial_t \bar{g}^{ij} &=& A^i_k \bar{g}^{kj}
+A^j_k \bar{g}^{ki} ,\\
 \partial_t \det \bar{g}_{ij} &=& -2 H \det
\bar{g}_{ij} ,
\\
 \partial_t h_{ij} &=& -\phi A_{ij} ,
\\
 \partial_t \det h_{ij} &=& - H \det h_{ij} ,
\\
 \partial_t K &=& HK ,
\\
\partial_t \phi &=& \frac{H}{n+2}\,\phi ,
\\
 \partial_t {g}_{ij} &=& -\frac{H}{n+2}\,
{g}_{ij}-A_{ij} ,
\\  \partial_t s &=& -\phi .
\end{eqnarray*}
\end{pro}

\begin{proof}
We interchange partial derivatives and use   equation (\ref{anf}) to
get
$$\partial_t
F_{,i}=\partial^2_{ti} F =\partial_i \xi = -A^k_i F_{,k} .$$ Note we
have also used the definition of affine curvature in equation
(\ref{afc}).

 Since $\partial_t\nu$ is a tangent vector,
\begin{eqnarray*}
\partial_t \nu &=& \langle \partial_t\nu , F_{,i}\rangle
\bar{g}^{ij}
F_{,j}\\
&= & -\langle  \nu , \partial^2_{ti} F \rangle
\bar{g}^{ij} F_{,j}\\
&= &  -\langle  \nu , -A^k_i F_{,k}\rangle \bar{g}^{ij}
F_{,j}\\
&= & 0.
\end{eqnarray*}

\begin{eqnarray*}
\partial_p \nu &= & \langle \partial_p \nu , F_{,i}\rangle
\bar{g}^{ij}
F_{,j}\\
&= & -\langle  \nu , \partial^2_{pi} F\rangle
\bar{g}^{ij} F_{,j}\\
&= &  -h_{pi} \bar{g}^{ij}F_{,j}.\\
\end{eqnarray*}

\begin{eqnarray*}
\partial_t \bar{g}_{ij} & = & \partial_t \langle \partial_i F,
\partial_j
F\rangle \\
&= &   \langle \partial^2_{ti} F , \partial_j F\rangle + \langle
\partial_i F , \partial^2_{tj}F\rangle\\
&= & \langle  -A^k_i \partial_k F , \partial_j F\rangle + \langle
\partial_i F
, -A^l_j \partial_l F\rangle\\
&= &  -A^k_i \bar{g}_{kj}  -A^k_j \bar{g}_{ki}.
\end{eqnarray*}

\begin{eqnarray*}
\partial_t \det\bar{g}_{ij}
&= & (\det\bar{g}_{lm})\bar{g}^{ij} \partial_t \bar{g}_{ij}\\
 &= &(\det\bar{g}_{lm})\bar{g}^{ij}(
-A^k_i \bar{g}_{kj}  -A^k_j \bar{g}_{ki}) \\
&= & -2 (\det\bar{g}_{lm}) H.
\end{eqnarray*}

\begin{eqnarray*}
\partial_t h_{ij} &= & \partial_t
\langle \partial^2_{ij} F , \nu \rangle \\
&= &   \langle \partial^3_{tij} F , \nu \rangle + \langle
\partial^2_{ij} F, \partial_t\nu \rangle\\
&= &\langle \partial^2_{ij} \xi , \nu \rangle\\
&= & \langle  \partial_i( -A^k_j \partial_k F) , \nu\rangle \\
&= &  -A^k_j h_{ik}.
\end{eqnarray*}

\begin{eqnarray*}
\partial_t \det h_{ij}
&= & (\det h_{lm}) h^{ij} \partial_t h_{ij}\\
 &= &(\det h_{lm}) h^{ij}(
  -A^k_j h_{ik}) \\
&= & - (\det h_{lm}) H
\end{eqnarray*}

Recall the formulas for the Gaussian curvature $K$, the affine
metric $g_{ij}$ and $\phi$:
$$K=
\frac{\det h_{ij}}{\det\bar{g}_{ij}}, \qquad  {g}_{ij}=
\frac{h_{ij}}{\phi}, \qquad \phi =K^{\frac{1}{n+2}}.$$

Thus, lowering the index on $A^k_i$ by the affine metric,
$$\partial_t h_{ij} =   -A^k_i h_{kj}  = - \phi h^{kl} A_{li} h_{kj}
=  -\phi A_{ij},
$$

\begin{eqnarray*}
\partial_t K &= & \partial_t \left(
\frac{\det h_{ij} }{\det \bar{g}_{ij} }\right )\\
&= & \frac{ (\partial_t \det h_{ij}) \det\bar{g}_{ij} -\det h_{ij}
(\partial_t\det \bar{g}_{ij}) }{(\det\bar{g}_{ij})^2}\\
&= & H K.
\end{eqnarray*}
and $$\partial_t \phi =  \frac{1 }{n+2} H \phi.
$$

Thus \begin{eqnarray*}\partial_t {g}_{ij}&=&\partial_t
\left(\frac{h_{ij}}{\phi}\right)\\
&=&(\partial_t h_{ij}) \left(\frac{1}{\phi}\right)-\frac{h_{ij}}{
\phi^2}\, \partial_t \phi\\
&=&(-\phi A_{ij}) \left(\frac{1}{\phi}\right)-\frac{h_{ij}}{
\phi^2}\left(\frac{1
}{n+2} H \phi\right)\\
&=& -\frac{H}{n+2} {g}_{ij}-A_{ij}.
 \end{eqnarray*}

$$\pt s = -\langle \pt  F , \nu
\rangle - \langle F,\pt\nu\rangle =-\langle \phi \nu , \nu \rangle -
0 = - \phi.$$
\end{proof}

\section{Evolution of the cubic form} \label{flow-C-sec}
We use the structure equation (\ref{s1}) to compute the evolution of
the cubic form.  First, we need to find the evolution of the affine
normal $\xi$ and of the Christoffel symbols.

\begin{pro}
Under the affine normal flow,
\begin{eqnarray*}
\pt \xi &=& -\frac1{n+2} g^{ij}H_{,i} \,F_{,j}
+ \frac{H}{n+2} \xi \\
&=& \frac1{n+2}\Delta \xi + \frac{2}{n+2}H \xi +
\frac{4}{n+2}A^m_iC^{ik}_{m}F_{,k}.\\
\end{eqnarray*}
\end{pro}

\begin{proof}
Recall $\xi = - g^{ki}(\ln \phi)_{,i} F_{,k} + \phi \nu$. First note
\begin{equation} \label{g-inverse-evolve}
\partial_t g^{iq} = -g^{i\ell}(\partial_t g_{\ell m}) g^{mq}
= - g^{i\ell} \left(-\frac{H}{n+2} g_{\ell m} - A_{\ell m} \right)
g^{mq}  = \frac{H}{n+2} g^{iq} + A^{iq}.
\end{equation}
Then compute using Proposition \ref{euclidean-evolve}
\begin{eqnarray*}
\pt \xi &=& \pt
\left(- g^{ki}(\ln \phi)_{,i} F_{,k} + \phi \nu\right) \\
&=& -(\partial_t g^{ki}) (\ln \phi)_{,i} F_{,k} - g^{ki}  (\pt\ln
\phi)_{,i}) F_{,k} - g^{ki} (\ln \phi)_{,i} (\partial_t F_{,k}) +
(\partial_t
\phi)\nu + 0 \\
&=& -\left(\frac{H}{n+2} g^{ki} + A^{ki}\right)(\ln \phi)_{,i}
F_{,k} - g^{ki} \left(\frac{H}{n+2}\right)_{,i} F_{,k} \\&&{}+
g^{ki} (\ln \phi)_{,i}
A^\ell_k F_{,\ell} + \frac{H}{n+2}\phi\nu \\
&=& -\frac1{n+2} g^{ij}H_{,i} F_{,j} + \frac{H}{n+2} \xi.
\end{eqnarray*}

From equation (\ref{F2cov}), we have $$\Delta \xi =
g^{ij}\xi_{,ij}=g^{ij}(-A^k_{i,j} F_{,k} - A_{ij}\xi - A^k_i
C^\ell_{kj} F_{,\ell})= -H \xi +g^{ij}(-A^k_{i,j} F_{,k}  - A^k_i
C^\ell_{kj} F_{,\ell}).$$ Now
$$g^{ij}A^k_{i,j}
F_{,k}=g^{ij}g^{kl}A_{il,j}F_{,k}=g^{ij}g^{kl}(A_{ij,l}+A^m_iC_{mlj}-A^m_l
C_{mij})F_{,k}=g^{kl}H_{,l} F_{,k}+A^m_iC^{ik}_{m}F_{,k}.$$ Hence
$$\Delta \xi =  -H \xi -g^{kl}H_{,l} F_{,k}
-2A^m_iC^{ik}_{m}F_{,k}$$ and
$$\heat \xi = \frac{2}{n+2}H \xi +
\frac{4}{n+2}A^m_iC^{ik}_{m}F_{,k}.$$
\end{proof}

We also compute
$$
\partial_t \Gamma^k_{ij} = \partial_t \sfrac12 g^{kl}(\partial_i
g_{j\ell} + \partial_j g_{i\ell} - \partial_\ell g_{ij})\\
$$
 Note $\partial_t \Gamma^k_{ij}$ is a tensor; therefore, we may
choose normal coordinates so that $\partial_k g_{ij} =
\Gamma^k_{ij} = 0$ at time $t=0$.  In these coordinates,
\begin{eqnarray*}
\partial_t \Gamma^k_{ij} &=& \sfrac12 g^{k\ell}\left[ \partial_i
\left( -\frac{H}{n+2} g_{j\ell} - A_{j\ell} \right) + \partial_j
\left(-\frac{H}{n+2} g_{i\ell} - A_{i\ell} \right)  \right. \\
 &&{}- \left.
\partial_\ell \left(-\frac{H}{n+2}g_{ij} - A_{ij} \right) \right]
\\
&=& -\frac1{2(n+2)}[(\partial_iH) \delta^k_j + (\partial_jH)
\delta^k_i - g^{k\ell} (\partial_\ell H) g_{ij}] \\
&&{} - \sfrac12 (\partial_i A^k_j + \partial_j A^k_i - g^{k\ell}
\partial_\ell A_{ij}) \\
&=& -\frac1{2(n+2)}(H_{,i} \delta^k_j + H_{,j} \delta^k_i -
g^{k\ell} H_{,\ell} g_{ij})  - \sfrac12 (A^k_{j,i} + A^k_{i,j} -
g^{k\ell} A_{ij,\ell})
\end{eqnarray*}

Now compute the evolution of $F_{,ij}$
\begin{eqnarray*}
\partial_t F_{,ij} &=& \partial_t\partial^2_{ij} F - (\partial_t
\Gamma^k_{ij})F_{,k} - \Gamma^k_{ij} \partial_t F_{,k} \\
&=& (\partial_t F)_{,ij} - (\partial_t \Gamma^k_{ij})F_{,k} \\
&=& \xi_{,ij}  + \frac1{2(n+2)}(H_{,i} \delta^k_j + H_{,j}
\delta^k_i - g^{k\ell} H_{,\ell} g_{ij}) F_{,k}  \\
&&{} + \sfrac12 (A^k_{j,i} + A^k_{i,j} - g^{k\ell} A_{ij,\ell})
F_{,k} \\
&=& -A^k_{i,j} F_{,k} - A_{ij}\xi - A^\ell_i C^k_{\ell j} F_{,k}
 + \sfrac12 (A^k_{j,i} + A^k_{i,j} - g^{k\ell} A_{ij,\ell})
F_{,k}\\
&& {} + \frac1{2(n+2)}(H_{,i} \delta^k_j + H_{,j} \delta^k_i -
g^{k\ell} H_{,\ell} g_{ij}) F_{,k}
\end{eqnarray*}
On the other hand,
\begin{eqnarray*}
\partial_t F_{,ij} &=& \partial_t (g_{ij} \xi + C^k_{ij} F_{,k})
\\
&=& \left(-\frac{H}{n+2} g_{ij} - A_{ij}\right) \xi + g_{ij}
\left(-\frac1{n+2}g^{k\ell} H_{,\ell} F_{,k} + \frac{H}{n+2} \xi
\right) \\
&&{} + (\partial_t C^k_{ij})F_{,k} - C^\ell_{ij} A_\ell^k F_{,k}
\end{eqnarray*}
Therefore,
\begin{eqnarray*}
\partial_t C^k_{ij} &=& -A^k_{i,j}  - A^\ell_i C^k_{\ell j}
+ \sfrac12 (A^k_{j,i} + A^k_{i,j} - g^{k\ell} A_{ij,\ell}) \\
&& {} + \frac1{2(n+2)}(H_{,i} \delta^k_j + H_{,j} \delta^k_i -
g^{k\ell} H_{,\ell} g_{ij})  \\
&&{} + C_{ij}^\ell A_\ell^k +\frac1{n+2} g_{ij}g^{k\ell}H_{,\ell}
\\
&=& -\sfrac12A_i^\ell C_{\ell j}^k - \sfrac12 A_j^\ell C^k_{\ell
i} -\sfrac12 g^{k\ell} A_{ij,\ell} \\
&& {} + \frac1{2(n+2)}(H_{,i} \delta^k_j + H_{,j} \delta^k_i -
g^{k\ell} H_{,\ell} g_{ij})  \\
&&{} + C_{ij}^\ell A_\ell^k +\frac1{n+2} g_{ij}g^{k\ell}H_{,\ell}
\end{eqnarray*}
The second line follows from the first by the Codazzi equation
(\ref{A-cod}) for $A^k_i$. Furthermore,
\begin{eqnarray*}
\partial_t C_{ijm} &=& \partial_t (g_{km} C^k_{ij}) \\
&=& \left(-\frac{H}{n+2}g_{km} - A_{km}\right)C^k_{ij}
-\sfrac12A_i^\ell C_{\ell jm} - \sfrac12 A_j^\ell C_{\ell
im} -\sfrac12 A_{ij,m} \\
&& {} + \frac1{2(n+2)}(H_{,i} g_{jm} + H_{,j} g_{im} -
 H_{,m} g_{ij})   + C_{ij}^\ell A_{\ell m} +\frac1{n+2}
 g_{ij}H_{,m}\\
&=& -\frac{H}{n+2} C_{ijm} + \frac1{2(n+2)}(H_{,i} g_{jm} + H_{,j}
g_{im} + H_{,m} g_{ij}) \\
&&{}-\sfrac12(A_{ij,m} - A_m^\ell C_{\ell ij}) -\sfrac12A_i^\ell
C_{\ell jm} - \sfrac12 A_j^\ell C_{\ell im} -\sfrac12A_m^\ell
C_{\ell ij}
\end{eqnarray*}
Note the first term in the last line is totally symmetric by the
Codazzi equation  (\ref{A-cod}) for $A^k_i$.

Now we compute the Laplacian of the cubic form.  We use apolarity
(\ref{apolarity}), the Codazzi equations (\ref{A-cod}) and
(\ref{C-cod}) for $A$ and $C$ respectively, and the curvature
equation (\ref{curv-eq}).
\begin{eqnarray*}
0 &=& g^{jk}C_{ijk,\ell m} \\
&=& g^{jk}(C_{i\ell k,jm} + \sfrac12 g_{ij}A_{k\ell,m} -\sfrac12
g_{i\ell}A_{jk,m} + \sfrac12g_{kj}A_{\ell i,m} -\sfrac12 g_{k\ell}
A_{ji,m})\\
&=& g^{jk} C_{i\ell j,km} +\sfrac12 A_{i\ell,m}
-\sfrac12g_{i\ell}H_{,m} + \sfrac12 n A_{\ell i,m} -\sfrac12
A_{\ell i, m} \\
&=& g^{jk}(C_{i\ell j,mk} - R_{kmip}C^p_{\ell j} - R_{km\ell p}
C^p_{ij} - R_{kmjp}C^p_{i\ell}) -\sfrac12g_{i\ell}H_{,m} +
\sfrac12 n A_{\ell i,m} \\
&=& g^{jk} C_{i\ell j,mk} -\sfrac12g_{i\ell}H_{,m} + \sfrac12 n
A_{\ell i,m} \\
&&{} - C_\ell^{pk}[\sfrac12(g_{ik}A_{mp}- g_{im}A_{kp} -
g_{pk}A_{mi} + g_{pm}A_{ki}) -C_{ik}^r C_{mrp} + C_{im}^r
C_{krp}]\\
&&{}- C_i^{pk}[\sfrac12(g_{\ell k}A_{mp}- g_{\ell m}A_{kp} -
g_{pk}A_{m\ell} + g_{pm}A_{k\ell}) -C_{\ell k}^r C_{mrp} + C_{\ell
m}^r C_{krp}] \\
&&{} - g^{jk}C_{i\ell p} [\sfrac12(g_{jk}A_{mp}- g_{jm}A_{kp} -
g_{pk}A_{mj} + g_{pm}A_{kj}) -C_{jk}^r C_{mrp} + C_{jm}^r
C_{jrp}]\\
&=& g^{jk} C_{ij\ell ,mk} -\sfrac12g_{i\ell}H_{,m} + \sfrac12 n
A_{\ell i,m} \\
&&{}+\sfrac12g_{im}A_k^pC^k_{p\ell} - \sfrac12A_i^jC_{m\ell j} + 2
C_{ik}^rC_{mr}^pC_{p\ell}^k - C_{im}^rC_{kr}^pC_{p\ell}^k +
\sfrac12 g_{\ell m} A_k^p C_{pi}^k \\
&&{} -\sfrac12 A_{k\ell}C^k_{mi} - C_{\ell m}^rC_{kr}^pC_{pi}^k -
\sfrac12 nA_m^pC_{pi\ell} -\sfrac12HC_{mi\ell} -
C_{jm}^rC_{rp}^jC_{i\ell}^p \\
&=& g^{jk}C_{i\ell m,jk} + \sfrac12 A_{m\ell,i} -\sfrac12 g_{im}
A^k_{\ell,k} + \sfrac12 A_{mi,\ell} - \sfrac12 g_{\ell m}
A^k_{i,k} -\sfrac12g_{i\ell}H_{,m} + \sfrac12 n A_{\ell i,m}\\
&&{}+\sfrac12g_{im}A_k^pC^k_{p\ell} - \sfrac12A_i^jC_{m\ell j} + 2
C_{ik}^rC_{mr}^pC_{p\ell}^k - C_{im}^rC_{kr}^pC_{p\ell}^k +
\sfrac12 g_{\ell m} A_k^p C_{pi}^k \\
&&{} -\sfrac12 A_{k\ell}C^k_{mi} - C_{\ell m}^rC_{kr}^pC_{pi}^k -
\sfrac12 nA_m^pC_{pi\ell} -\sfrac12HC_{mi\ell} -
C_{jm}^rC_{rp}^jC_{i\ell}^p
\end{eqnarray*}
Now the Codazzi equation (\ref{A-cod}) for $A^k_i$  and the
apolarity condition (\ref{apolarity}) imply $$A^k_{i,k} = A^k_{k,i}
+ A_k^\ell C_{\ell i}^k - A_i^\ell C^k_{\ell k} = H_{,i} +
A^p_kC_{pi}^k.$$
 Apply this identity and the Codazzi equation (\ref{A-cod})
for $A^k_i$ to the first two occurrences of the covariant
derivatives of $A$ to find
\begin{eqnarray*}
\Delta C_{i\ell m} &=& g^{jk}C_{i\ell m,jk} = \sfrac12
g_{im}H_{,\ell} + \sfrac12 g_{\ell m} H_{,i} + \sfrac12 g_{i\ell}
H_{,m}
+\sfrac12(n+2)(A_m^kC_{ki\ell} - A_{\ell i,m}) \\
&&{}- 2C_{ik}^rC_{mr}^pC_{p\ell}^k + C_{im}^r C_{kr}^p C_{p\ell}^k
+ C_{\ell m}^r C_{kr}^p C_{pi}^k + C_{jm}^rC_{rp}^jC_{i\ell}^p +
\sfrac12 H C_{mi\ell}
\end{eqnarray*}
Together with the evolution equation of $C$, compute
\begin{eqnarray*}
\partial_t C_{ijk} &=& \frac1{n+2}\Delta C_{ijk}
-\frac{3H}{2(n+2)}C_{ijk} + \frac2{n+2}
C_{i\ell}^mC_{km}^pC_{pj}^\ell \\
&&{} -\frac1{n+2}(C_{ik}^m C_{\ell m}^p C_{pj}^\ell + C_{jk}^m
C_{\ell m}^p C_{pi}^\ell + C_{\ell k}^m C_{mp}^\ell C_{ij}^p) \\
&&{}-\frac12(A_i^\ell C_{\ell jk} + A_j^\ell C_{\ell ik} +
A_k^\ell C_{\ell ij})
\end{eqnarray*}

To compute $\partial_t |C|^2$, use (\ref{g-inverse-evolve}) to show
\begin{eqnarray*}
\partial_t|C|^2 &=& \partial_t (C_{ijk} g^{iq}g^{jr}g^{ks}
C_{qrs}) \\
&=& 3 C_{ijk}(\partial_t g^{iq}) C_q^{rs} + 2
(\partial_t C_{ijk}) C^{ijk} \\
&=& \frac{3H}{n+2} |C|^2 + 3C_{ik}^rA^{iq}C_{qr}^k + \frac2{n+2}
\Delta C_{ijk} C^{ijk} - \frac{3H}{n+2} |C|^2 \\
&&{}+ \frac4{n+2} C_{i\ell}^mC_{km}^p C_{pj}^\ell C^{ijk} -
\frac6{n+2} C_{ik}^m C_{\ell m}^p C^\ell_{pj} C^{ijk} - 3A^{i\ell}
C^j_{\ell k} C^k_{ij} \\
&=&   \frac2{n+2} \Delta C_{ijk} C^{ijk} + \frac4{n+2}
C_{i\ell}^mC_{km}^p C_{pj}^\ell C^{ijk} - \frac6{n+2} |P|^2
\end{eqnarray*}
for $P_{ij} = C^k_{i\ell} C^\ell_{jk}$. Finally compute $\Delta
|C|^2 = 2 \Delta C_{ijk} C^{ijk} + 2 |\nabla C|^2$ to find
 $$\partial_t |C|^2 = \frac1{n+2} \Delta |C|^2 - \frac2{n+2}
 |\nabla
 C|^2 + \frac4{n+2}
C_{i\ell}^mC_{km}^p C_{pj}^\ell C^{ijk} - \frac6{n+2} |P|^2.$$

Now if $\mathcal Y_{ijkl} = C_{ij}^mC_{klm} - C_{ik}^mC_{jlm}$, we
find $$0\le\frac12|\mathcal Y|^2 = |P|^2 - C_{i\ell}^mC_{km}^p
C_{pj}^\ell C^{ijk},$$ and so
$$ \partial_t |C|^2 \le \frac1{n+2}\Delta|C|^2 - \frac2{n+2}|P|^2
\le \frac1{n+2}\Delta|C|^2 - \frac{2}{n(n+2)}|C|^4,$$ since
$|C|^2=P^i_i$ and thus Cauchy-Schwartz applied to the eigenvalues of
$P$ implies $|P|^2\ge \frac1n |C|^4$. We note this estimate of
Andrews \cite{andrews96} is a parabolic version of an estimate of
Calabi \cite{calabi72} on the cubic form on affine spheres, and is
related to Calabi's earlier interior $C^3$ estimates of solutions to
the Monge-Amp\'ere equation \cite{calabi58}.

The maximum principle implies the following estimate for $|C|^2$
then: If $\mathcal L$ is any compact smooth strictly convex
hypersurface evolving as $\mathcal L(t)$ under the affine normal
flow, then
$$ \sup_{\mathcal L(t)} |C|^2 \le \frac 1{(\sup_{\mathcal
L(0)}|C|^2)^{-1} + \frac2{n(n+2)}\,t}.$$ Thus we get the following
bound independent of initial data:
\begin{pro}[Andrews \cite{andrews96}] \label{c2-decay}
Let $\mathcal L$ be any compact smooth strictly convex hypersurface
evolving under the affine normal flow. Then
$$ \sup_{\mathcal L(t)} |C|^2 \le \frac{n(n+2)}{2t}.$$
\end{pro}

\section{Evolution of the affine curvature} \label{flow-A-sec}

In this section, we treat the evolution of the affine curvature
$A^i_k$, as computed by Andrews \cite{andrews96}, and also the
evolution of the affine conormal vector $U$. At each point, $U$ is
defined by
 \begin{equation} \label{define-conormal} \langle U,\xi \rangle = 1, \qquad
\langle U,F_{,i}\rangle = 0, \quad i=1,\dots,n.
\end{equation}
(It should be clear that in this case, we are using the Euclidean
inner product $\langle\cdot,\cdot\rangle$ only for notational
convenience. As its name suggests, the conormal vector $U$ is more
naturally a vector in the dual space to $\R^{n+1}$, not a vector in
$\R^{n+1}$ itself.)

\begin{pro} Under the affine normal flow,
\begin{eqnarray*}
\pt U &=& -\frac H{n+2}\,U = \frac1{n+2}\Delta U, \\
 \pt A^k_i &=& A_i^j A_j^k
+\frac1{n+2}H_{,\ell i} \,g^{\ell k} + \frac1{n+2} H_{,\ell}
\,g^{\ell
j} C_{ji}^k + \frac{H}{n+2} A^k_i,\\
\pt A_{ij} &=& \frac1{n+2} H_{,ij} +
\frac1{n+2} H_{,\ell}\,g^{\ell k} C_{ijk}, \\
\pt H &=& \frac1{n+2}\Delta H + |A|^2 +
\frac1{n+2} H^2.\\
 \pt A_{ij} &=&\frac1{n+2}\Delta  A_{ij} -\frac{1}{n+2} \Big(2A^{pk}C_{pmk}C^{m}_{ij}+ 2A^{ml}
C_{mij,l}+A^p_iC_{pml}C^{ml}_{j}\\
&&{}+C^{lp}_{i} C_{lpm}A^m_j- 2A^p_lC_{pmi}C^{ml}_{j}  -
g_{ij}A^k_{m}A^m_k  + nA^m_iA_{mj}\Big)
\end{eqnarray*}
\end{pro}

\begin{proof}
Compute $\pt U$ by differentiating its defining equation
(\ref{define-conormal}):
 \begin{eqnarray*}
 \langle \partial_t U,\xi\rangle &=& - \langle
U,\partial_t \xi\rangle = - \left\langle
-\frac1{n+2}g^{ij}H_{,i}F_{,j} + \frac H{n+2}\,\xi \right\rangle = -
\frac H{n+2}. \\
\langle \partial_t U, F_{,i} \rangle &=& - \langle U,
\partial_t F_{,i} \rangle = -\langle U, - A_i^kF_{,k}\rangle = 0, \\
\partial_t U &=& - \frac H{n+2}\,U.
\end{eqnarray*}
Similarly, covariantly differentiate in space to find
\begin{eqnarray*}
\langle U_{,i}, \xi \rangle &=& - \langle U,\xi_{,i} \rangle = -
\langle U, -A_i^k F_{,k} \rangle = 0, \\
\langle U_{,i}, F_{,j} \rangle &=& - \langle U, F_{,ij}\rangle = -
\langle U, g_{ij}\xi + C_{ij}^k F_{,k} \rangle = -g_{ij}, \\
\langle U_{,ij},\xi\rangle &=& - \langle U_{,i}, \xi_{,j} \rangle =
-\langle U_{,i}, -A_j^l F_{,l} \rangle = -g_{il}A^l_j = -A_{ij},\\
\langle U_{,ij},F_{,k}\rangle &=& - \langle U_{,i},F_{,kj} \rangle =
- \langle U_{,i},g_{kj}\xi + C^l_{kj} F_{,l} \rangle =
g_{il}C^l_{kj} = C_{ijk}.
\end{eqnarray*}
Now for $\Delta U = g^{ij} U_{,ij}$, we have $\langle g^{ij}
U_{,ij}, \xi \rangle = -H$, $\langle g^{ij} U_{,ij}, F_{,k} \rangle
= g^{ij}C_{ijk}=0$ by the apolarity of the cubic form. So $\Delta U
= -HU$ and
$$\partial_t U = \frac1{n+2}\, \Delta U.$$

 To compute $\partial_t A^i_k$, we use the defining equation for $A$:
$\xi_{,i} = -A^k_i F_{,k}$. Take $\partial_t$ to find
$$ \left(-\frac1{n+2}H_{,\ell}g^{\ell k} F_{,k} + \frac{H}{n+2} \xi
\right)_{,i} = \partial_t \xi_{,i} = -\partial_t (A^k_i F_k)  =
-(\partial_t A^k_i) F_{,k} + A^k_iA^j_k F_{,j}.$$
 So we have
\begin{eqnarray*}
 (\partial_tA^k_i) F_{,k} &=& A^k_iA^j_k F_{,j} -
 \left(-\frac1{n+2}H_{,\ell}g^{\ell k} F_{,k} + \frac{H}{n+2} \xi
\right)_{,i} \\
&=&A^k_iA^j_k F_{,j} -\frac1{n+2}[-H_{,\ell i} g^{\ell k} F_{,k} -
H_{,\ell} g^{\ell k} F_{,ki}
\\&&{} + H_{,i} \xi - H
A_i^k F_{,k}] \\
&=& \left(A^j_iA_j^k  + \frac1{n+2} H_{,\ell i} g^{\ell k} +
\frac1{n+2} H_{,\ell} g^{\ell j} C^k_{ji} +\frac{H}{n+2}
A_i^k\right) F_{,k}.
\end{eqnarray*}
Here we have used the structure equation (\ref{s1}).
\begin{eqnarray} \nonumber
\partial_t A_{im} &=& \partial_t(g_{km}A_i^k) \\
\nonumber
&=& g_{km}\left(A^j_iA_j^k  + \frac1{n+2} H_{,\ell i} g^{\ell k} +
\frac1{n+2} H_{,\ell} g^{\ell j} C^k_{ji} +\frac{H}{n+2}
A_i^k\right) \\
 \nonumber && {}+ A_i^k \left(-\frac{H}{n+2} g_{km} -
A_{km}\right) \\
 \label{pt-Aij} &=&\frac1{n+2} H_{,mi} + \frac1{n+2}
H_{,\ell} C^\ell_{im}.
\end{eqnarray}
Finally,
\begin{eqnarray*}
\partial_t H &=& \partial_t A^i_i \\
&=&A_i^j A_j^i +\frac1{n+2}H_{,\ell i} \,g^{\ell i} + \frac1{n+2}
H_{,\ell} \,g^{\ell j} C_{ji}^i + \frac{H}{n+2} A^i_i \\
&=& \frac1{n+2} \Delta H + |A|^2 + \frac{H^2}{n+2}
\end{eqnarray*}
by the apolarity condition $C_{ji}^i=0$.
\end{proof}

\begin{eqnarray*}
\Delta  A_{ij} &=& g^{kl} A_{ij,kl}= g^{kl} A_{ji,kl} \\
&=&  g^{kl} (A^m_k C_{mij} + A_{jk,i} -A^m_iC_{mjk})_l \\
&=& g^{kl} (A_{mk,l} C^{m}_{ij} + A^m_k C_{mij,l} + A_{jk,il}-
A_{mi,l} C^m_{jk}-A^m_i C_{mjk,l}) \\
\end{eqnarray*}

\begin{eqnarray*}
 A_{jk,il}&=&  A_{jk,li}+R^m_{il j} A_{mk}+ R^m_{i lk} A_{jm}\\
 &=& A_{kj,li}
-\Big[ \sfrac12(- g_{lj} A^{m}_i +g_{ij} A^{m}_l +\delta^{m}_l
A_{ij}-\delta^{m}_i A_{lj}) +
C^p_{lj}C^{m}_{pi}-C^p_{ij}C^{m}_{pl}\Big]A_{mk}\\
&-&  \Big[ \sfrac12(- g_{lk} A^{m}_i +g_{ik} A^{m}_l +\delta^{m}_l
A_{ik}-\delta^{m}_i A_{lk}) +
C^p_{lk}C^{m}_{pi}-C^p_{ik}C^{m}_{pl}\Big]A_{mj}\\
\end{eqnarray*}

\begin{eqnarray*} A_{jk,li} &=& A_{kl,ji}+(A^m_lC_{mjk}- A^m_jC_{mlk})_i\\
&=& A_{kl,ji}+(A_{ml,i}C^m_{jk} + A^m_lC_{mjk,i}- A_{mj,i}C^m_{lk}
-A^m_jC_{mlk,i} ) \\
\end{eqnarray*}

\begin{eqnarray*}
 g^{kl}A_{jk,il}
 &=&
H_{,ij}+A_{ml,i}C^{ml}_{j} + A^m_lC^l_{mj,i} -\sfrac12
g_{ij}A^m_{l}A^l_m+ \sfrac{n}2A^m_iA_{mj}\\
&-& C^p_{lj} C^m_{pi}A^l_{m} +C^p_{ij} C^m_{pl} A^l_{m}+C^p_{ik}
C^k_{pm}A^m_j
\end{eqnarray*}

\begin{eqnarray*}
A^m_lC^l_{mj,i} &=& A^{mk}C_{mkj,i} \\
&=& A^{mk} C_{mjk,i}\\
&=& A^{mk} [C_{mji,k}+\sfrac{1}{2}
(g_{mk}A_{ji}+g_{jk}A_{mi}-g_{mi}A_{jk}-g_{ji}A_{mk})]\\
&=& A^{mk} C_{mij,k}+\sfrac{1}{2}
HA_{ij}-\sfrac{1}{2}g_{ij}A_{mk}A^{mk}\\
\end{eqnarray*}

\begin{eqnarray*}
 g^{kl}A_{jk,il}
 &=&
H_{,ij}+A_{ml,i}C^{ml}_{j}+A^{mk} C_{mij,k}+\sfrac{1}{2}
HA_{ij}-\sfrac{1}{2}g_{ij}A_{mk}A^{mk}
 -\sfrac12
g_{ij}A^m_{l}A^l_m\\&&{}+ \sfrac{n}2A^m_iA_{mj} - C^p_{lj}
C^m_{pi}A^l_{m} +C^p_{ij} C^m_{pl} A^l_{m}+C^p_{ik} C^k_{pm}A^m_j
\end{eqnarray*}

\begin{eqnarray*}
\Delta A_{ij} &=& g^{kl} (A_{mk,l} C^{m}_{ij} + A^m_k C_{mij,l} +
A_{jk,il}-
A^m_{i,l} C_{mjk}-A^m_i C_{mjk,l}) \\
&=& g^{kl} (A_{mk,l} C^{m}_{ij} + A^m_k C_{mij,l} + A_{jk,il}-
A^m_{i,l} C_{mjk}-A^m_i C_{mjk,l}) \\
\end{eqnarray*}

By the Codazzi equations (\ref{A-cod}) and (\ref{C-cod}),
\begin{eqnarray*}
g^{kl}A_{mk,l} C^{m}_{ij}&=&H_{,m}C^{m}_{ij}+
A^{pk}C_{pmk}C^{m}_{ij},\\
 g^{kl}A^m_{l,i} C_{mjk}&=&g^{kl}A_{mi,l}C^m_{jk}+A^p_iC_{mpl}
C^{ml}_{j}-A^p_lC_{mpi} C^{ml}_{j}, \\
g^{kl}A^m_iC_{mjk,l}&=& \frac{1}{2} (HA_{ij} +A^m_iA_{mj} -
A^m_iA_{mj}-nA^m_iA_{mj})= \frac{1}{2} (HA_{ij} -nA^m_iA_{mj}).
\end{eqnarray*}

Hence
\begin{eqnarray*}
\Delta A_{ij} &=& g^{kl} (A_{mk,l} C^{m}_{ij} + A^m_k C_{mij,l}
+ A_{jk,il}- A^m_{i,l} C_{mjk}-A^m_i C_{mjk,l} )\\
&=& H_{,m}C^{m}_{ij}+ A^{pk}C_{pmk}C^{m}_{ij}+ A^{ml}
C_{mij,l}+H_{,ij}+A_{ml,i}C^{ml}_{j}+A^{mk}
C_{mij,k}\\
&&{}+\sfrac{1}{2} HA_{ij}-\sfrac{1}{2}g_{ij}A_{mk}A^{mk}
 -\sfrac12
g_{ij}A^m_{l}A^l_m+ \sfrac{n}2A^m_iA_{mj} - C^p_{lj} C^m_{pi}A^l_{m}
\\&&{}+C^p_{ij} C^m_{pl} A^l_{m}+C^p_{ik}
C^k_{pm}A^m_j
-A_{mi,l} C^{ml}_{j}- \sfrac{1}{2} (HA_{ij} -nA^m_iA_{mj})\\
&=& H_{,ij}+H_{,m}C^{m}_{ij}+ 2A^{pk}C_{pmk}C^{m}_{ij}+ 2A^{ml}
C_{mij,l}+A^p_iC_{pml}C^{ml}_{j}+ A^m_jC_{lpm}C^{lp}_{i}\\
&&{}- 2A^{pl}g^{mn}C_{pmi}C_{lnj}  - g_{ij}A^k_{m}A^m_k
+nA^m_iA_{mj}
\end{eqnarray*}

By the evolution equation (\ref{pt-Aij}) of $A_{ij}$,
\begin{eqnarray*}
 \pt A_{ij} &=&\frac1{n+2}\Delta  A_{ij} -\frac{1}{n+2} \Big(2A^{pk}C_{pmk}C^{m}_{ij}+ 2A^{ml}
C_{mij,l}+A^p_iC_{pml}C^{ml}_{j}\\&&{}+C^{lp}_{i} C_{lpm}A^m_j -
2A^p_lC_{pmi}C^{ml}_{j}  - g_{ij}A^k_{m}A^m_k  + nA^m_iA_{mj}\Big).
\end{eqnarray*}



\section{The support function} \label{support-function-section}
In this section, we recall some standard facts about the support
function of a convex body in $\R^{n+1}$, derive the equation
satisfied by the support function under the affine normal flow, and
use convexity to prove local $C^0$ and $C^1$ estimates for support
functions of a family of smooth bounded convex domains exhausting a
general convex domain.

Below we will consider the following situation: Let $\mathcal K =
\bigcup_{i=1}^\infty \mathcal K^i$ be a convex domain in
$\R^{n+1}$ exhausted by bounded convex domains $\mathcal K^i$. Our
initial hypersurface $\mathcal L = \partial \mathcal K$ will then
be considered as a limit of the more regular hypersurfaces
$\mathcal L^i = \partial \mathcal K^i$. Let $\mathcal L^i(t)$ and
$\mathcal L(t)$ denote the affine normal flow with initial
hypersurface $\mathcal L^i$ and $\mathcal L$ respectively.

Then, for an initial convex hypersurface $\mathcal L =
\partial \mathcal K = \partial \left(\bigcup_{i=1}^\infty \mathcal
K^i \right),$ we want local uniform estimates of the affine normal
flow $\mathcal L^i(t)$ as $\mathcal L^i(t)\to \mathcal L(t)$. In
this section, we recall some standard facts about the support
function and use convexity to prove $C^0$ and $C^1$ estimates
locally in $\mathcal D^\circ(s_{\mathcal K})$.

Recall for $\mathcal L = \partial \mathcal K$, the support function
is defined for $Y\in\R^{n+1}$ by $$s(Y) = \sup_{x\in \mathcal K}
\langle x,Y \rangle.$$

Here are some important properties of the support function (see
Rockafellar \cite{rockafellar}). First of all, recall equation
(\ref{deriv-s}) that in the case $\mathcal L$ is smooth and
strictly convex, the total derivative of the support function $ds
= F$ the embedding.  In our case, $\mathcal L$ may not be smooth
and strictly convex; but we may still recover the convex domain
$\mathcal K$ from the support function. Take the \emph{Legendre
transform} of $s$ : For $x\in \R^{n+1}$, let
$$\delta(x) = \sup_{Y\in\R^{n+1}} \langle x,Y\rangle-s(Y).$$
Then $\delta$ is the \emph{indicator function} of the closed convex
set $\bar{\mathcal K}$.  In other words, $$\delta(x) =
\left\{ \begin{array}{c@{\,\,\mbox{ for }}l} 0 & x\in \bar {\mathcal K} \\
+\infty & x\notin \bar {\mathcal K} . \end{array} \right.$$

Let $\mathcal D(s) = s^{-1}(-\infty,+\infty)\subset \R^{n+1}$ be the
\emph{domain} of the support function $s$, and let $\mathcal
D^\circ(s)$ denote the interior of the domain.  The support function
of a convex domain $\mathcal K$ is always a convex,
lower-semicontinuous function $s\!:\R^{n+1}\to (-\infty,+\infty]$ of
homogeneity one.   Moreover, any convex lower-semicontinuous
function $s\!:\R^{n+1}\to (-\infty,+\infty]$ of homogeneity one is
the support function of a closed convex set so long as $s$ is not
identically $+\infty$. The support function, since it is convex, is
continuous on $\mathcal D^\circ(s)$ but may not be continuous on all
of  $\mathcal D(s)$.

The following lemma follows from the description above of the
Legendre transform of the support function:
\begin{lem} \label{compare-support}
If $Q_1$ and $Q_2$ are closed convex subsets of $\R^{n+1}$, then
$Q_1 \subset Q_2$ if and only if the support functions $s_{Q_1}\le
s_{Q_2}$ on all $\R^{n+1}$.
\end{lem}

All of our estimates will be uniform on compact subsets of $\mathcal
D^\circ(s)\times (0,T]$ for some positive time $T$.  So we need the
following lemma to start

\begin{lem}
If $\mathcal K$ is a convex domain in $\R^{n+1}$ which contains no
lines, then for the support function $s_{\mathcal K}$, $\mathcal
D^\circ(s_{\mathcal K})\neq \emptyset$.
\end{lem}

\begin{proof}
We prove the lemma by contradiction.  If $\mathcal
D^\circ(s_{\mathcal K})= \emptyset$, then then since $\mathcal
D(s_{\mathcal K})$ is a convex collection of rays, $\mathcal
D(s_{\mathcal K})$ must be contained in a hyperplane $\mathcal
H=\{Y:\langle Y,v\rangle =0\}$. Since $s_{\mathcal K}|_{\mathcal H}$
is a convex function of homogeneity one on $\mathcal H$, there is a
linear function $\langle Y,w\rangle$ which is $\le s$ on $\mathcal
H$. Now consider the line $L=\{w+\tau v : \tau\in\R\}$, whose
support function is
$$s_L(Y)=\left\{
\begin{array}{c@{\mbox{ for }}c} +\infty & \langle Y,v \rangle \neq
0 \\ \langle Y,w\rangle & \langle Y,v \rangle = 0. \end{array}
\right.$$
 By construction, $s_L\le s_{\mathcal K}$ on $\R^{n+1}$, and so
$L\subset \bar {\mathcal K}$ by Lemma \ref{compare-support}. The
convex hull of $L$ and any open ball in $\mathcal K$ then contains
another line contained in the open set $\mathcal K$, and this
provides a contradiction.
\end{proof}

\begin{pro} \label{support-exhaust-limit}
Let $$\mathcal K = \bigcup_{i=1}^\infty \mathcal K^i$$ be convex
bodies so that $\mathcal K^i \subset \mathcal K^{i+1}$. Then the
support functions $s = s_{\mathcal K}$, $s_i = s_{\mathcal K^i}$
satisfy $s_{i+1} \ge s_i$ and $s_i\to s$ everywhere, and the
convergence is uniform on compact subsets of $\mathcal D^\circ (s)$.
If, in addition, each $\mathcal K^i$ is bounded with smooth,
strictly convex boundary, then the $C^1$ norm of $s_i$ is uniformly
bounded on each compact subset of $\mathcal D^\circ(s)$.
\end{pro}

\begin{proof}
First of all, it is clear from the definition of $s$ that
$s_{i+1}\ge s_i$, and $s(Y) = \lim_{i\to\infty} s_i(Y)$ for all
$Y\in \R_{n+1}$:
$$
s(Y) = \sup_{x\in \mathcal K} \langle x,Y\rangle  =
\sup_{x\in\bigcup\mathcal K^i} \langle x,Y \rangle  = \sup_i
\sup_{x\in \mathcal K^i} \langle x,Y\rangle  = \sup_i s_{\mathcal
K^i}(Y)  = \lim_{i\to\infty} s_{i}(Y) $$
 since $\{s_i(Y)\}$ is an
increasing sequence for all $Y$.

Let $C\subset \mathcal D^\circ$ be a compact subset.  Choose a
compact $C'\subset \mathcal D^\circ$ which contains a neighborhood
of $C$.  Note that on all of $\mathcal D^\circ$, for all $i$, $$s_1
\le s_i \le s.$$  Thus for $Y\in\partial C$, we have
$$|ds_i(Y)| \le \frac{\max_{\partial C'} \left|s\right| -
\min_{\partial C} \left|s_1\right|}{\mbox{dist }(\partial
C',\partial C)}.$$ (Proof: For every direction $v$, consider $s_i$
restricted to the line $L$ through $Y$ with direction $v$. Then the
directional derivative of $s_i$ at $Y$ is bounded above by the slope
of the secant line of the graph of $s_i$  through $Y$ and a point in
$L\cap \partial C'$.)

Since each $s_i$ is convex, the same estimate is true on all of $C$.
Therefore, the $C^1$ norm of all the $s_i$ is bounded on $C$, and
since we have pointwise convergence, Ascoli-Arzel\'a implies uniform
convergence of $s_i\to s$ on $C$.
\end{proof}

Now we recall the standard formulas for the support function of a
domain with smooth and strictly convex boundary, in particular,
relating it to the Gauss curvature.  We also derive the parabolic
Monge-Amp\`ere equation the support function satisfies under the
affine normal flow.

Recall above that $\partial_t s = -\phi= - K^{\frac1{n+2}}$.  We now
derive some standard formulas relating the Gauss curvature $K$ to
the support function $s$.

Recall $s(Y)$ is a convex function on $\R^{n+1}$ which is
homogeneous of degree one.  Let $F(x)$ denote a local embedding of a
smooth, strictly convex hypersurface $\mathcal L =
\partial \mathcal K$.  Then at any $F(x)\in \mathcal L$ at
which $s(Y)=\langle F(x),Y\rangle$, $Y$ is perpendicular to the
tangent space $T_{F(x)}\mathcal L$.  By restricting to $Y$ on the
unit sphere $\sph^n$ in $\R^{n+1}$, we have a natural
parametrization of $\mathcal L$, which is given by the inverse of
the Gauss map $-\nu$. For $F(x)\in \mathcal L$, let $Y=-\nu(x)$ be
the outward normal.  Then since $\mathcal L$ is strictly convex,
$x\mapsto Y$ is a local diffeomorphism for $Y\in\sph^n$, and we can
consider $F=F(Y)$ for $Y\in\sph^n$.  We extend $F$ to be homogeneous
of order zero:
$$F\!: \R^{n+1}\setminus \{0\} \to \R^{n+1}, \qquad F(Y) =
F\left(\frac{Y}{|Y|} \right).$$ Then $s(Y) = -\langle F,\nu \rangle$
and thus
\begin{equation} \label{support-fy}
s(Y) = \langle F,Y\rangle
\end{equation}
for all $Y\in\R^{n+1}\setminus\{0\}$.

It is useful to consider the support function restricted to an
affine hyperplane of distance 1 to the origin in $\R^{n+1}$.  We may
choose coordinates so that $$Y=(y,-1)=(y^1,\dots,y^n,-1).$$  By
projecting from this hyperplane to $\sph^n$, we still have a local
parametrization of our hypersurface $\mathcal L$, and
(\ref{support-fy}) still holds.  Now differentiate
(\ref{support-fy}) to find for $i=1,\dots,n$
$$
\frac{\partial s}{\partial y^i} = \left\langle \frac{\partial
F}{\partial y^i}, Y \right\rangle + F^i = F^i$$
 since $Y$ is normal to $\mathcal L$.  Moreover, we use Euler's
 formula $$\sum_{i=1}^{n+1} y^i \frac{\partial F}{\partial y^i} = 0$$ to
 show \begin{eqnarray*}
 \frac{\partial F}{\partial y^{n+1}} &=& -\frac1{y^{n+1}}\sum_{i=1}^n y^i
 \frac{\partial F}{\partial y^i} = \sum_{i=1}^n y^i
 \frac{\partial F}{\partial y^i}, \\
\frac{\partial s}{\partial y^{n+1}} &=&
 \left\langle \frac{\partial
F}{\partial y^{n+1}}, Y \right\rangle + F^{n+1} \\
 &=& \left\langle
\sum_{i=1}^n y^i \frac{\partial F}{\partial
y^i},Y \right\rangle + F^{n+1} \\
&=& F^{n+1}
 \end{eqnarray*}
since $Y$ is normal to the image of $F$. Thus at any
$Y\in\R^{n+1}\setminus \{0\}$, the total derivative
 \begin{equation} \label{deriv-s} ds =
(F^1,\dots,F^{n+1}) = F. \end{equation}

Now differentiate $\langle \frac{\partial F}{\partial y^i}, Y\rangle
= 0$ to find for $i,j=1,\dots,n$,
\begin{eqnarray*}
 0 &=& \frac{\partial}{\partial y^j} \left\langle
\frac{\partial F}
{\partial y^i}, Y \right\rangle \\
 &=& \left\langle \frac{\partial^2 F}{\partial y^i
\partial y^j}, Y \right\rangle + \frac{\partial F^j}{\partial
y^i},\\
\frac{\partial^2 s}{\partial y^i
\partial y^j} = \frac{\partial F^j}{\partial y^i} &=& -\left\langle
\frac{\partial^2 F}{\partial y^i
\partial y^j}, Y \right\rangle \\
 &=& \left\langle \frac{\partial^2
F}{\partial y^i
\partial y^j}, \nu|Y| \right\rangle \\
&=& \sqrt{1+|y|^2} \left\langle \frac{\partial^2 F}{\partial y^i
\partial y^j}, \nu\right\rangle  \\
&=& h_{ij}\sqrt{1+|y|^2}.
\end{eqnarray*}

Moreover, we compute for $i,j=1,\dots,n$
\begin{eqnarray*}
\bar g_{ij} &=& \left\langle \frac{\partial F}{\partial y^i},
\frac{\partial F}{\partial y^j} \right\rangle \\
&=& \sum_{k=1}^{n+1} \frac{\partial F^k}{\partial y^i}
\frac{\partial F^k}{\partial y^j} \\
&=& \frac{\partial F^{n+1}}{\partial y^i} \frac{\partial
F^{n+1}}{\partial y^j} +  \sum_{k=1}^n \frac{\partial^2 s}{\partial
y^i
\partial y^k}
\frac{\partial^2s}{\partial y^j \partial y^k} \\
&=& \frac{\partial F^i}{\partial y^{n+1}} \frac{\partial
F^j}{\partial y^{n+1}} + \sum_{k=1}^{n} \frac{\partial^2 s}{\partial
y^i \partial y^k}
\frac{\partial^2s}{\partial y^j \partial y^k} \\
&=& \left(\sum_{k=1}^n \frac{\partial F^i}{\partial y^k} y^k \right)
\left(\sum_{l=1}^n \frac{\partial F^j}{\partial y^l} y^l \right) +
\sum_{k=1}^{n} \frac{\partial^2 s}{\partial y^i \partial y^k}
\frac{\partial^2s}{\partial y^j \partial y^k} \\
&=& \sum_{k,l=1}^n \frac{\partial^2s}{\partial y^i \partial y^k}
(y^ky^l + \delta^{kl}) \frac{\partial^2s}{\partial y^j \partial
y^l}, \\
\det{\bar g_{ij}} &=& \det \left(\frac{\partial^2s}{\partial y^i
\partial y^k} \right) \det (y^ky^l + \delta^{kl}) \det \left(
\frac{\partial^2s}{\partial y^j \partial y^l}\right) \\
&=& (1+|y|^2)\det \left(\frac{\partial^2 s}{\partial y^i
\partial y^j} \right)^2.
\end{eqnarray*}
So the Gaussian curvature
 \begin{eqnarray*}
K &=& \frac{\det h_{ij}}{\det \bar g_{ij}} =
(1+|y|^2)^{-\frac{n+2}2} \left( \det \frac{\partial^2 s}{\partial
y^i \partial y^j} \right)^{-1}, \\
\phi &=& K^{\frac1{n+2}} = (1+|y|^2)^{-\frac12} \left( \det
\frac{\partial^2 s}{\partial y^i \partial y^j}
\right)^{-\frac1{n+2}}.
\end{eqnarray*}

In order to address the evolution of $s$, we note \emph{a priori}
that there are two natural parametrizations $F$ of our hypersurface.
First, the affine normal flow defines a particular parametrization
at time $t>0$ given an initial parametrization at time $t=0$. On the
other hand, for any hypersurface $F(y,t)$, there is a natural
parametrization in terms of the inverse of the Gauss map $-\nu$.
These two parametrizations are compatible in the following sense:

\begin{pro}
Given a hypersurface $\mathcal L \subset \R^{n+1}$ parametrized by
the inverse of its Gauss map $F\!: \sph^n\to \mathcal L$, under the
affine normal flow, $F(y,t)$ is still a parametrization by the
inverse of the Gauss map.
\end{pro}
\begin{proof}
The two parametrizations are related by the Gauss map $-\nu$. Under
the affine normal flow,  $\nu$ satisfies $\partial_t\nu=0$ by
Proposition \ref{euclidean-evolve}.
\end{proof}

Thus if we assume the initial parametrization is via the inverse of
the Gauss map, the formulas developed in this section are still
valid under the affine normal flow (and in any case the two
parametrization merely differ by a diffeomorphism).

Denote by $s(y)$ $$s(y) =  s(y^1,\dots,y^n,-1) = \sqrt{1+|y|^2}
\,\,s\left(\frac Y{|Y|} \right)$$ for $Y/|Y|\in \sph^n$.  Thus we
find under the affine normal flow
$$
\partial_t s(y) =
\sqrt{1+|y|^2}\,\,
\partial_t s\left(\frac Y{|Y|} \right) = -\phi\sqrt{1+|y|^2} = -
\left(\det \frac{\partial^2s}{\partial y^i
\partial y^j} \right)^{-\frac1{n+2}},
$$
where we have used $\partial_t s = -\phi$ from Proposition
\ref{euclidean-evolve}.

We record this as
\begin{pro}
For any smooth solution to the affine normal flow, the support
function $s(y)$ as defined above satisfies
 \begin{equation} \label{s-y-flow}
 \partial_t s(y) = -
\left(\det \frac{\partial^2s}{\partial y^i
\partial y^j} \right)^{-\frac1{n+2}}.
\end{equation}
\end{pro}

\section{The flow} \label{flow-def-sec}

There is no question about the definition of affine normal flow
beginning at a smooth strictly convex compact hypersurface in
$\R^{n+1}$ (this is true for any convex compact hypersurface by
Andrews \cite{andrews00}).  It is convenient to define the affine
normal flow for an open convex domain in $\R^{n+1}$ by performing
affine normal flow on the boundary of the domain.  In this way we
let $\Psi_t \mathcal J = \mathcal J(t)$ denote the affine normal
flow of $\mathcal J$ a bounded domain with smooth strictly convex
boundary in $\R^{n+1}$. For $t$ larger than the extinction time,
define $\Psi_t \mathcal J = \mathcal J(t) = \emptyset$.

Note to pass from a convex embedded hypersurface $\mathcal L$ to a
domain $\mathcal J$ with $\mathcal L=\partial \mathcal J$, set
$\mathcal J$ to be the interior of the convex hull $(\hat {\mathcal
L})^\circ$.

Consider an open convex region $\mathcal K\subset \R^{n+1}$ which
contains no lines. Then the boundary $\partial \mathcal K$ is a
properly embedded convex hypersurface in $\R^{n+1}$.  We define the
affine normal flow on the hypersurface $\partial \mathcal K$ by its
action on the interior of its convex hull $\mathcal K$. Now we
define the affine normal flow on the hypersurface $\partial \mathcal
K$ and on the region $\mathcal K$ by
 \begin{equation} \label{flow-def}
  \mathcal K(t) = \bigcup_{\mathcal J \subset \mathcal K } \mathcal J(t),
 \end{equation}
where each $\mathcal J$ in (\ref{flow-def}) is a bounded domain with
smooth strictly convex boundary.

\begin{lem} \label{compact}
If $\mathcal L$ is a compact convex hypersurface in $\R^{n+1}$, then
our definition of the affine normal flow $\mathcal L(t)$ corresponds
with the usual one.
\end{lem}

\begin{proof}
If $\mathcal L$ is strictly convex and smooth, then this follows
at once from the maximum principle.  Otherwise, Andrews
\cite{andrews00} shows that there is a viscosity solution
$\tilde{\mathcal L}(t)$ to the affine normal flow which is unique
provided that the Hausdorff distance from $\tilde {\mathcal L}(t)$
to $\mathcal L$ goes to zero as $t\to0$. Moreover $\tilde
{\mathcal L}(t)$ is smooth and strictly convex for positive $t$
less than the extinction time.

Our definition $\mathcal L(t)$ is clearly a viscosity solution, and
the Hausdorff convergence property is satisfied by Lemma \ref{cont}
below. Therefore, Andrews's uniqueness result implies $\tilde
{\mathcal L}(t) = \mathcal L(t)$, and so our definition coincides
with the standard one in the compact case.
\end{proof}

\begin{rem}
We recall (see e.g.\ \cite{andrews00}) that a \emph{viscosity
solution} to a hypersurface flow problem is a family of
hypersurfaces $\mathcal L(t)$ with initial condition $\mathcal
L(0)=\mathcal L$ so that: 1) If $\mathcal J$ is a smooth
hypersurface contained in $\mathcal L$, then the evolving
hypersurface $\mathcal J(t)$ is contained in $\mathcal L(t)$ for all
$t\in[0,T]$, and 2) If $\mathcal J$ is a smooth hypersurface
containing $\mathcal L$, then the evolving hypersurface $\mathcal
J(t)$ contains $\mathcal L(t)$ for all $t\in[0,T]$.  In short, a
viscosity solution $\mathcal L(t)$ is one for which the maximum
principle always works, even if $\mathcal L(t)$ does not have $C^2$
regularity.
\end{rem}

The following proposition depends on estimates proved by Ben
Andrews in the case of compact hypersurfaces \cite{andrews00}.
Below, we prove local versions of the estimates needed.
\begin{pro} \label{exhaust-flow}
Let $\mathcal K^i$ and $\mathcal K$ be open convex bodies containing
no lines so that
$$\mathcal K=\bigcup_{i=1}^\infty \mathcal K^i,
\qquad \mathcal K^i\subset \mathcal K^{i+1}.$$
Then for all $t>0$,
$$\mathcal K(t) = \bigcup_{i=1}^\infty \mathcal K^i(t).$$
\end{pro}

\begin{proof}
Let $\mathcal J\subset \mathcal K$ be a bounded domain with smooth,
strictly convex boundary.  Then
$$\mathcal J = \bigcup_{i=1}^\infty
\mathcal J^i, \qquad \mathcal J^i = \mathcal J \cap \mathcal K^i.$$
Then we claim that
\begin{equation} \label{compact-exhaust}
\mathcal J(t) = \bigcup_{i=1}^\infty \mathcal J^i(t)
\end{equation}

To prove the claim (\ref{compact-exhaust}), we recall estimates of
Andrews \cite[Section 8]{andrews00} for compact convex
hypersurfaces (we prove local versions of these estimates below).

By exhausting $\mathcal J^i = \bigcup_{j=1}^\infty \mathcal I^i_j$
by nested domains $\mathcal I^i_j$ with smooth, strictly convex
boundary, the affine normal flow $\mathcal J^i(t)$ is defined as a
limit as $j\to\infty$ of the affine normal flow $\mathcal I^i_j(t)$.
The support functions $s_{\mathcal I^i_j} \to s_{\mathcal J^i}$
uniformly on compact subsets of $\R^{n+1}$ as $j\to\infty$.  The
resulting $C^0$ estimates automatically entail parabolic $C^{2,1}$
estimates for positive $t$ (see below), and then Krylov's theory
implies parabolic $C^{2+\alpha,1+\frac\alpha2}$ estimates. These
estimates ensure that the limit $s_{\mathcal J^i}(t)$ of the
$s_{\mathcal I^i_j}(t)$ exists and is smooth for $t>0$.  Andrews
shows this solution is unique by applying barriers and the maximum
principle.

The key point is that $C^0$ estimates on the support function of
convex bounded regions imply local parabolic
$C^{2+\alpha,1+\frac\alpha2}$ estimates of the affine normal flow
for all times $t>0$. Since $\mathcal J = \bigcup_{i=1}^\infty
\mathcal J^i$, we have that $s_{\mathcal J^i} \to s_{\mathcal J}$
locally in $C^0$.  Therefore, under the affine normal flow,
$s_{\mathcal J^i}(t)$ converges to a limit $s(t)$ locally in
parabolic $C^{2+\alpha,1+\frac\alpha2}$ for $t>0$.  Since
$\mathcal J = \bigcup_{i=1}^\infty \mathcal J^i$, the limit $s(t)$
converges uniformly on convex sets to $s_{\mathcal J}(0)$ as
$t\to0$.  Andrews's uniqueness argument then shows that
$s(t)=s_{\mathcal J}(t)$ and the claim (\ref{compact-exhaust}) is
proved.

Now use (\ref{compact-exhaust}) to compute
$$
\mathcal K(t) = \bigcup_{\mathcal J\subset \mathcal K} \mathcal J(t)
= \bigcup_{\mathcal J\subset \mathcal K} \left( \bigcup_{i=1}^\infty
\mathcal J^i(t) \right)  \subset  \bigcup_{\mathcal J\subset
\mathcal K} \left( \bigcup_{i=1}^\infty \mathcal K^i(t) \right)  =
\bigcup_{i=1}^\infty \mathcal K^i(t).
$$
($\mathcal J$ of course represents bounded domains with smooth,
strictly convex boundaries.) On the other hand,
$$
\mathcal K(t) = \bigcup_{\mathcal J\subset\mathcal K} \mathcal J(t)
= \bigcup_{\mathcal J \subset \bigcup_{i=1}^\infty \mathcal K^i}
\mathcal J(t)  \supset \bigcup_{i=1}^\infty \left(\bigcup_{\mathcal
J\subset \mathcal K^i} \mathcal J(t) \right)  = \bigcup_{i=1}^\infty
\mathcal K^i(t).
$$
This completes the proof of  Proposition
\ref{exhaust-flow}.
\end{proof}

The following corollary ensures convexity
\begin{cor}
$\mathcal K(t)$ is convex for all $t>0$ before the extinction time.
\end{cor}
\begin{proof}
Choose each $\mathcal K^i$ in the previous proposition to be a
bounded domain with strictly convex smooth boundary. Then $\mathcal
K(t)$ is an increasing union of convex sets.
\end{proof}

We verify that our definition satisfies the semigroup property:
\begin{lem} $\Psi_t \Psi_s \mathcal K = \Psi_{t+s} \mathcal K$.
\end{lem}
\begin{proof}
We work in terms of the support functions.  Let $s_{\mathcal
K}(Y,t)$ denote the support function of the domain $\Psi_t\mathcal
K$.  We claim
\begin{equation} \label{support-semigroup}
s_{\mathcal K}(Y,t+s) = s_{\Psi_s \mathcal K} (Y,t).
\end{equation}

To prove the claim, write $\mathcal K = \bigcup_{i=1}^\infty
\mathcal K^i$, where each $\mathcal K^i$ is a bounded domain with
smooth strictly convex boundary and $\mathcal K^i \subset \mathcal
K^{i+1}$ for all $i$. Since the semigroup property holds for each
$\mathcal K^i$, we have
$$\Psi_{t+s}\mathcal K^i = \Psi_t\Psi_s\mathcal K^i \quad
\Longrightarrow \quad s_{\mathcal K^i}(Y,t+s) = s_{\Psi_s\mathcal
K^i}(Y,t)$$ for all $t,s>0$ and $Y\in\R^{n+1}$.  Now let
$i\to\infty$. Propositions \ref{support-exhaust-limit} and
\ref{exhaust-flow} then prove the claim (\ref{support-semigroup}).

The lemma follows from (\ref{support-semigroup}) because any open
convex domain can be recovered from its support function by taking
the Legendre transform \cite{rockafellar}.
\end{proof}

We also have a lemma on the continuity of the flow:
\begin{lem} \label{cont}
For any $\tau\ge0$, and $\mathcal K$ a convex body in $\R^{n+1}$
containing no lines,
$$ \mathcal K(\tau) = \bigcup_{t>\tau}\mathcal K(t).$$
\end{lem}

\begin{proof}
By the semigroup property, we may assume $\tau=0$.  Consider any
point $p\in\mathcal K$. Since $\mathcal K$ is open, there is a small
ball around $p$ contained in $\mathcal K$.  This ball acts as a
barrier under the affine normal flow, and $p\in \mathcal K(t)$ for
$t>0$ the extinction time of the affine normal flow of this ball.
\end{proof}

By means of outer barriers, we show our definition actually
corresponds to the usual definition of affine normal flow for a
smooth, strictly convex hypersurface. Let Aff$(n+1)$ denote the
special affine group $\mathbf{SL}(n+1)\ltimes \R^{n+1}$.

\begin{pro} \label{def-works} Let
$\mathcal L\subset \R^{n+1}$ be a properly embedded convex
hypersurface which contains no lines. Assume in a neighborhood of a
point $p\in \mathcal L$ that $\mathcal L$ is $C^2$ and strictly
convex. Then
$$ \frac{\partial \mathcal L}{\partial t} (p) =
\xi_p \,\,\mbox{mod} \,\, T_p(\mathcal L(t)).$$
\end{pro}

\begin{rem}
This proposition should also follow from the estimates proved below
(what is still needed in addition is a local version of Andrews's
speed bound).
\end{rem}

\begin{proof}
Note that of course the  derivative $\frac{\partial \mathcal
L}{\partial t}(p)$ is defined only when $\mathcal L(t)$ is locally
parametrized. This parametrization defines the derivative, but
different parametrizations may cause it to vary by an element of
the the tangent space $T_p(\mathcal L(t))$.

Since $\Psi_t$ is a semigroup, we may assume $t=0$.

To proceed with the proof, we need a lemma on choosing nice
coordinates.

\begin{lem}
Let $p\in \mathcal L\subset \R^{n+1}$, and let $\mathcal L$ be a
$C^2$ strictly convex hypersurface near $p$. Then there is an
element $\Phi\in\mathbf{SL}(n+1,\R)\ltimes \R^{n+1}$ so that
$p\mapsto 0$ and the image locally is
 \begin{equation} \label{K-graph} \Phi (\mathcal L) =
 \left\{ x^{n+1} = \frac\gamma 2 |x|^2 +
o(|x|^2) \right\}
\end{equation}
 for $x=(x^1,\cdots,x^n)$, $\gamma>0$.
\end{lem}
\begin{proof}
This amounts to using Aff$(n+1)$ to choose coordinates.  Use a
rotation to set the inward-pointing normal to be $e_{n+1}$, and
translate so that $p$ is at the origin. We can still move the
tangent plane $\{x^{n+1}=0\}$ by an action of $\mathbf{SL}(n,\R)$.
Since $\mathcal L$ is strictly convex, we have
$$ \mathcal L = \left\{x^{n+1} = \sum_{i,j}a_{ij}x^i x^j + o(|x|^2)
\right\}$$ for $(a_{ij})$ a positive definite symmetric matrix. Use
the action of $\mathbf{SL}(n,\R)$ to send the ellipsoid
$a_{ij}x^ix^j \le C$ to a sphere of the same volume.  This amounts
to setting $a_{ij}= \frac\gamma2 \delta_{ij}$ for a positive
constant $\gamma$.
\end{proof}

Locally we write  $\mathcal L(t) = \{x^{n+1} = f(t,x)\}$ for
$f(0,x)=f(x)$ given in (\ref{K-graph}).  Modulo a tangential
piece, the affine normal to $\mathcal L$ at 0 is $\xi = (\det
f_{ij})^\frac{1}{n+2}e_{n+1}$ (see e.g.\ Nomizu-Sasaki, page 48).
Thus we want to show that
$$ \frac{\partial f}{\partial t} (0) = (\det
f_{ij})^\frac{1}{n+2}.$$
 In other words, we want to show
 \begin{equation} \label{squeeze}
  \limsup_{t\to0^+} \frac{f(t,0)-f(0,0)}t \le
 (\det f_{ij})^\frac{1}{n+2} \le
\liminf_{t\to0^+} \frac{f(t,0)-f(0,0)}t.
 \end{equation}

The left-hand inequality in (\ref{squeeze}) is equivalent to
showing that for each smooth strictly convex hypersurface
$\mathcal H\subset \hat {\mathcal L}$, $p\in \mathcal H$, the
affine normal of $\mathcal H$ at $p$ is $\le (\det
f_{ij})^\frac{1}{n+2}$. This is true by the definition. Consider
$\mathcal H$ a compact, strictly convex, $C^2$ hypersurface so
that $\mathcal H \subset \hat {\mathcal L}$ and $\mathcal H$
coincides with $\mathcal L$ in a neighborhood of the point $p$.
Then the definition gives $\hat {\mathcal L}(t) \supset \mathcal
H(t)$ for all small positive $t$. Therefore, the left-hand
inequality in (\ref{squeeze}) is proved.

To show the right-hand inequality in (\ref{squeeze}), we find
specific hyperboloid barriers whose $e_{n+1}$ component of the
affine normal at $p$ approaches $(\det f_{ij})^\frac{1}{n+2}$.

So choose $\epsilon>0$. Then consider hyperboloids of the form
$$\{x_{n+1} = G(x) = \sqrt{\alpha|x|^2 + \beta} - \sqrt\beta\}$$
for $\alpha,\beta>0$. Then compute the Hessian matrix $G_{ij}(0) =
\beta^{-\frac12}\alpha \delta_{ij}$. Fix $\alpha$ and $\beta$ so
that $\beta^{-\frac12}\alpha = \gamma - \epsilon$ so that
$$ G(x)=G_\beta(x) = \sqrt{\sqrt\beta(\gamma-\epsilon) |x|^2 + \beta}
- \sqrt\beta.$$ Then for $x$ in a small ball $B$ near 0,
$G_1(x)\le f(x)$ by (\ref{K-graph}). Now since $\mathcal L$ is
convex, $\mathcal L\setminus B$ lies above the graph of a function
$c|x|$ for $c$ a positive constant. $G_\beta(x)\to 0$ as $\beta\to
0^+$, and moreover $\partial G_\beta / \partial \beta \ge0$. So we
may choose $\beta$ close to zero so that $G_\beta(x)\le f(x)$ on
$B$ and also $G_\beta(x)\le c|x|$. Therefore, $\mathcal L$ lies
above the graph of $G_\beta(x)$, and the graph of $G_\beta$ is a
barrier to all compact hypersurfaces contained in $\mathcal L$.
The $e_{n+1}$ component of the affine normal of the graph of
$G_\beta$ at 0 is given by
$$\left(\det \frac{\partial^2G_\beta}{\partial x^i
\partial x^j}(0)\right)^{\frac1{n+2}} =
(\gamma-\epsilon)^{\frac{n}{n+2}}.$$ Therefore, by our definition of
affine normal flow, we have
$$\liminf_{t\to0^+} \frac{f(t,0)-f(0,0)}t \ge
(\gamma-\epsilon)^\frac{n}{n+2}$$ for all $\epsilon>0$. Now let
$\epsilon\to 0$ and Proposition \ref{def-works} is proved.
\end{proof}

An important and easy consequence of our definition is the following

\begin{thm} [Maximum principle at infinity]
Consider two convex domains $\mathcal K^1\subset \mathcal K^2$. Then
for all positive $t$, $\mathcal K^1(t) \subset \mathcal K^2(t)$.
\end{thm}

\begin{rem}
There are other natural flows for which such a global maximum
principal fails.  For example, there is an example in Ecker
\cite{ecker97} in which two spacelike soliton solutions to the mean
curvature flow in Minkowski space cross at infinity in finite time.
\end{rem}

\begin{thm}[Long Time Existence] \label{long-time-ex}
Let $\mathcal K$ be an unbounded convex domain in $\R^{n+1}$ which
contains no lines. Then for all $t>0$, $\mathcal K(t)\neq
\emptyset$.
\end{thm}

\begin{proof}
It is well known that such a $\mathcal K$ contains an infinite
half-cylinder in $\R^{n+1}$. Therefore, $\mathcal K$ contains
ellipsoids of unbounded volume, which in turn have unbounded
extinction times (the ellipsoids are equivalent under the action of
Aff$(n+1)$ to spheres of unbounded volume; these have unbounded
extinction times, as in Example \ref{sph-example}). The maximum
principle then completes the proof.
\end{proof}

\begin{pro}
If $\mathcal K$ is a convex domain in $\R^{n+1}$ which contains a
line, then the affine normal flow leaves $\mathcal K$ unchanged.
\end{pro}

\begin{proof}
Let $p\in \mathcal K$ and recall $\mathcal K$ is open.  Then
$\mathcal K$ contains a round cylinder centered at $p$.  This
cylinder contains ellipsoids of arbitrarily large volume centered at
$p$, which act as barriers to the affine normal flow.  These
barriers ensure that $p$ is always in $\mathcal K(t)$.  Thus
$\mathcal K(t)=\mathcal K$ for all $t\ge0$.
\end{proof}

\section{Soliton solutions} \label{sol-ex-sec}
It is well-known that solitons of the affine normal flow are the
convex properly embedded affine spheres---see Proposition
\ref{soliton-affine-sphere} below. These were classified by Cheng
and Yau \cite{cheng-yau86}.

\begin{pro} \label{soliton-affine-sphere}
Under the affine normal flow, an expanding soliton is a hyperbolic
affine sphere, a translating soliton is a parabolic affine sphere,
and a shrinking soliton is an elliptic affine sphere.
\end{pro}

\begin{proof}
This is a simple local calculation. $F=F(x)$ is a local embedding of
an expanding soliton which expands away from a central point $P$ at
a given point in time if and only if $$ \partial_t F = \xi =
\lambda(F-P) + Z^iF_{,i}$$ for $\lambda$ a positive constant and
$Z^iF_{,i}$ a tangent vector field.

The equiaffine condition of the affine normal (as in Proposition
\ref{aff-normal-def}), states that $\xi_{,j}$ is contained in the
tangent space.  Thus we compute
$$\xi_{,j} = \lambda F_{,j} + Z^iF_{,ij} + Z^i_{,j} F_{,i} =
(\lambda\delta^i_j + Z^i_{,j})F_{,i} + Z^i(g_{ij}\xi +
C^k_{ij}F_{,k}).$$ By comparing both sides of the equation in the
span of $\xi$, we find $Z^ig_{ij}=0$, and so the tangential piece
$Z^i=0$. Therefore, $\xi = \lambda(F-P)$, which is the equation
for a hyperbolic affine sphere centered at $P$. The cases of
shrinking and translating solitons are essentially the same.
\end{proof}

So shrinking solitons are elliptic affine spheres, and the only
properly embedded examples are ellipsoids \cite{cheng-yau86}, which
are images under Aff$(n+1)$ of the Euclidean spheres discussed
above. Since they are compact, Lemma \ref{compact} shows our
definition corresponds with the classical one.  We record the
example of the round sphere.

\begin{exa} \label{sph-example}
For a sphere of radius $r$ in $\R^{n+1}$, the affine normal $\xi$ is
$r^{-\frac n{n+2}}$ times the unit inward normal vector.  So then,
the affine normal flow $\partial_t F = \xi$ becomes the ODE $dr/dt =
-r^{-\frac{n}{n+2}}$, and so the radius at time $t$ satisfies
$$ r(t) = \left( r_0^{\frac{2n+2}{n+2}} - \frac{2n+2}{n+2}
\,t\right) ^{\frac{n+2}{2n+2}}.$$ The extinction time of a sphere
with initial radius $r_0$ is $$ \frac{n+2}{2n+2}
r_0^{\frac{2n+2}{n+2}}.$$ Note that if the initial radius (or the
initial enclosed volume) of a family of spheres tends to $\infty$,
then the extinction time goes to $\infty$.
\end{exa}

Translating solitons are parabolic affine spheres, and the only
properly embedded examples are elliptic paraboloids.  Expanding
solitons are hyperbolic affine spheres, and for every convex cone in
$\R^{n+1}$ containing no lines, there is a homothetic family of
hyperbolic affine spheres asymptotic to the cone (for the standard
round cone, these are simply hyperboloids).  In the next few
examples, we verify that our definition of affine normal flow leads
to the correct behavior for these solitons.

\begin{exa} We consider the paraboloid $\mathcal L=\{x^{n+1} = |x|^2\}$.  Our
affine normal flow $\Psi_t$ is invariant under the action of ${\rm
Aff}(n+1)$.  Consider a point $P = (\tilde x^1,\dots,\tilde
x^{n+1})$ on the paraboloid.  Then the following map in ${\rm
Aff}(n+1)$ preserves the paraboloid:
$$ (x^1,\dots,x^n,x^{n+1}) \mapsto \left(x^1+\tilde x^1, \dots, x^n
+ \tilde x^n, x^{n+1} + \tilde x^{n+1} + 2 \sum_{i=1}^n  x^i \tilde
x^i \right).$$ This map sends the origin to $P$, and also sends
$$(0,\dots,0, c) \mapsto P + (0,\dots,0,c).$$  Since $\Psi_t$ is invariant
under such transformations, each $\Psi_t \mathcal L=\mathcal L(t)$
must be a paraboloid $x^{n+1} = |x|^2 + c(t)$, and since $\Psi_t$ is
the usual affine normal flow pointwise, $\mathcal L(t)$ is the
standard translating soliton.
\end{exa}

\begin{exa} \label{weak-hyp-as}
Let $\mathcal C$ be a convex cone in $\R^{n+1}$ which contains no
lines and has the origin as vertex. Then $\mathcal C$ is invariant
under scaling by positive constants.  Of course, such homothetic
scalings are not in ${\rm Aff}(n+1)$ in general, but we can still
use the scaling properties of the usual affine normal flow to
determine how $\Psi_t \mathcal C = \mathcal C(t)$ scales in time.

It is straightforward to check that for each compact convex region
$\mathcal K$ and scale parameter $\lambda>0$,
$$
 \Psi_t (\lambda
\mathcal K) = \lambda \Psi_{t\lambda^{-(2n+2)/(n+2)}} \mathcal K.
$$
Because
$$ \mathcal K \subset \mathcal C \iff \lambda \mathcal K \subset \mathcal C,$$ and by
our definition of $\Psi_t$, then
$$ \Psi_t(\mathcal C) = \Psi_t (\lambda \mathcal C) =
\lambda \Psi_{t\lambda^{-(2n+2)/(n+2)}} \mathcal C$$
 for all $\lambda>0$.  Thus for all $t>0$, we have
 $$ \Psi_t \mathcal C = t^{\frac{n+2}{2n+2}} \Psi_1 \mathcal C.$$
Thus the hypersurface evolving from any such cone $\mathcal C$
scales as a hyperbolic affine sphere under the affine normal flow.
Indeed, this expanding soliton is (homothetic scalings of)
Cheng-Yau's hyperbolic affine sphere.  The local regularity results
below will prove that this viscosity solution is the same as
Cheng-Yau's smooth solution.
\end{exa}

If the cone $\mathcal C$ is homogeneous, then we can use the full
affine symmetry group to conclude much more.  Below we analyze the
affine normal flow of Calabi's example \cite{calabi72}. The case
of the hyperboloid is similar (the symmetry group being the
Lorentz group in this case).

\begin{exa} \label{calabi-ex}
Let $\mathcal C = \mathcal C(0)$ be the boundary of the first
orthant. Then since $\mathcal C$ is invariant under multiplying all
the coordinates by positive scalars, and the flow is invariant under
Aff$(n+1)$, $\mathcal C(t)$ is invariant under the group
$$G=\{(\lambda_1,\dots,\lambda_{n+1}) : \lambda_i>0, \prod \lambda_i
= 1 \}$$ acting by
$$(x^1,\dots,x^{n+1}) \mapsto (\lambda_1x^1,\dots, \lambda_{n+1}
x^{n+1}).$$
 At time $\epsilon$, Calabi's example does move (consider a
 hyperboloid containing the first orthant with vertex at the
origin).  By group invariance, $\mathcal C(t)$ must be of the form
$$\left\{\prod_i x^i = {\rm const.},\, x^i>0\right\},$$ which is an
orbit of $G$. Proposition \ref{def-works} shows that at time $t>0$,
our flow is the affine normal flow pointwise.  The radial graph of
$\Psi_t \mathcal C$ must solve an ODE in $t$, and so it must be the
standard solution
 $$\mathcal C(t) = \{(x^1,\dots,x^{n+1})\in\R^{n+1} :
 x^i\ge0, \prod_i x^i  = c_n t^{\frac{n+2}2}\}$$
 for $c_n = (n+1)^\frac12 (\frac2{n+2})^{\frac{n+2}2}.$
\end{exa}

\section{Andrews's Speed Estimate} \label{andrews-est-sec}

The following proposition is essentially found in Andrews
\cite{andrews00}, following work of Tso \cite{tso85}, although the
statement of the proposition in \cite{andrews00} is slightly
incorrect.

\begin{pro} \label{andrews-speed}
Let $s$ be the support function of a smooth strictly convex compact
hypersurface evolving under affine normal flow. If $s(Y,t)\ge r>0$
for all $Y\in\sph^n$ and $t\in [0,T]$, then
$$|\partial_t s| \le \left(C + C' t^{-\frac n{2n+2}}\right)s$$
on $\sph^n\times[0,T]$, where $C$ and $C'$ are constants only
depending on $r$ and $n$.
\end{pro}

\begin{proof}
Consider the function $$q = \frac{-\partial_t s}{s - r/2}.$$
 We apply
the maximum principle to $\log q = \log |\partial_ts| -
\log(s-r/2)$.  In particular, at a fixed time $t\in[0,T]$, consider
a point $Y\in\sph^n$ at which $q$ attains its maximum. By changing
coordinates, we may assume that this point $Y=(0,\dots,0,-1)$ is the
south pole.  Then, as in Tso, consider the coordinates
$y=(y^1,\dots,y^n)$ for $s$ restricted to the hyperplane
$\{(y^1,\dots,y^n,-1)\}$.  At $y=0$, we have for $i=1,\dots,n$
\begin{equation}
(\log q)_i =0 \quad \Longleftrightarrow \quad \frac{s_{ti}}{s_t} =
\frac{s_i}{s-r/2} \label{grad-zero}
\end{equation}
The condition for $(\log q)|_{\sph^n}$ to have a maximum at the
south pole is
\begin{equation}
 (\log q)_{ij} + (\log q)_{n+1}
\delta_{ij}\le 0 \label{neg-def}
\end{equation} as a symmetric matrix.
Here we use subscripts to denote ordinary differentiation
$f_i=\partial_{y^i}f$ and $f_t = \partial_t f$.

To compute the second term in (\ref{neg-def}), use Euler's
identities
$$\sum_{i=1}^{n+1} y^i s_{ti} = s_t, \qquad \sum_{i=1}^{n+1}
y^i s_i = s$$ at the point $Y=(0,\dots,0,-1)$ to conclude $s_{tn+1}
= -s_t$, $s_{n+1} = -s$,  and
$$ (\log q)_{n+1} =
\frac{r/2}{s-r/2}.
$$
For the first term in (\ref{neg-def}), compute
\begin{eqnarray*}
(\log q)_{ij} &=& \frac{s_{tij}}{s_t} - \frac{s_{ti}s_{tj}}{s_t^2} -
\frac{s_{ij}}{s-r/2} + \frac{s_is_j}{(s-r/2)^2} \\
&=& \frac{s_{tij}}{s_t}  - \frac{s_{ij}}{s-r/2}
\end{eqnarray*}
at $y=0$ by (\ref{grad-zero}).  Thus (\ref{neg-def}) becomes at
$y=0$
\begin{equation}
\label{neg-def2} \frac{r/2}{s-r/2}\delta_{ij} + \frac{s_{tij}}{s_t}
- \frac{s_{ij}}{s-r/2} \le 0.
\end{equation}

Now, we compute using the flow equation (\ref{s-y-flow})
\begin{eqnarray*}
(\log q)_t &=& \partial_t \log |\partial_t s| - \partial_t
\log(s-r/2) \\
&=& -\frac1{n+2}\,\partial_t \log\det (s_{ij}) - \frac{s_t}{s-r/2} \\
&=& -\frac1{n+2} \,s^{ij} s_{tij} - \frac{s_t}{s-r/2}
\end{eqnarray*}
for $s^{ij}$ the inverse matrix of $s_{ij}$. Then (\ref{neg-def2})
implies that
\begin{eqnarray*}
(\log q)_t &\le&  \frac{r/2}{n+2} \cdot \frac{s_t}{s-r/2}
\delta_{ij} s^{ij} - \frac{2n}{n+2} \cdot \frac{s_t}{s-r/2} \\
&=& -\frac{r/2}{n+2} \,q \,\delta_{ij}s^{ij} + \frac{2n}{n+2}\, q, \\
q_t &\le& -\frac{r/2}{n+2} \,q^2\, \delta_{ij}s^{ij} +
\frac{2n}{n+2} \,q^2.
\end{eqnarray*}
Now if we let $\mu_i$ be the eigenvalues of $s^{ij}$, or
equivalently the reciprocals of the eigenvalues of $s_{ij}$, then we
see
$$ |s_t| = (\det s_{ij})^{-\frac1{n+2}} = \left(\prod_{i=1}^n \mu_i
\right)^{\frac1{n+2}} \le \left(\frac1n \sum_{i=1}^n \mu_i
\right)^{\frac n{n+2}} = \left(\frac1n \,\delta_{ij}s^{ij}
\right)^{\frac n{n+2}}$$ by the arithmetic-geometric mean
inequality.  Therefore,
$$ \delta_{ij}s^{ij} \ge n |s_t|^{\frac{n+2}n} = n q^{\frac{n+2}n}
(s-r/2)^{\frac{n+2}n} \ge n q^{\frac{n+2}n} (r/2)^{\frac{n+2}n}$$
since $s\ge r$.  And so finally, at $y=0$, and thus at any maximum
point of $q|_{\sph^n}$,
\begin{equation} \label{qt-at-max}
q_t \le -\frac{n(r/2)^{\frac{2n+2}n}}{n+2}\, q^{\frac{3n+2}n} +
\frac{2n}{n+2} \, q^2.
\end{equation}

Now define $Q(t) = \max_{Y\in \sph^n} q(Y,t)$.  Then
(\ref{qt-at-max}) implies that
$$ Q_t \le -Q^2\left(c_n r^{\frac{2n+2}n} Q^{\frac{n+2}n}
-c_n'\right)$$ for constants $c_n,c_n'$ depending only on $n$.
Therefore,
 \begin{equation} \label{Q-bound}
 Q\le \max \left\{c_n r^{-\frac{2n+2}{n+2}}, c_n' r^{-1}
t^{-\frac n{2n+2}}\right\}
 \end{equation}
for $c_n,c_n'$ new constants depending only on $n$. The result
easily follows.

\begin{rem} $Q$ may not be differentiable as a function of $t$, but the
above estimate (\ref{Q-bound}) still holds---see e.g.\ Hamilton
\cite[Section 3]{hamilton86}.
\end{rem}

\end{proof}

\section{Guti\'errez-Huang's Second Derivative Estimate}
\label{gh-est-sec}

In this section, we follow Guti\'errez-Huang
\cite{gutierrez-huang98} to find an upper bound of the Hessian of
solutions to  the equation satisfied by the support function under
the affine normal flow
$$\partial_t s = - \left(\det \frac{\partial^2 s}{\partial y^i
\partial y^j} \right)^{-\frac1{n+2}}.$$
This is a parabolic version of an estimate of Pogorelov for
elliptic Monge-Amp\`ere equations.  We will treat the slightly
more general case
\begin{equation} \label{gen-eq}
 \partial_t u = -\rho(y) \left(\det
\frac{\partial^2 u}{\partial y^i
\partial y^j} \right)^{-\alpha},
\end{equation}
for $\rho(y)$ a smooth positive function on $\R^n$ and $\alpha$ a
positive constant. Guti\'errez and Huang considered the case
$\rho(y)=\alpha=1$.  The reason we introduce $\rho(y)$ is that the
evolution of the support function of a hypersurface by a power of
the Gauss curvature involves a term $\rho(y)$ which is a power of
$1+|y|^2$. The calculations are essentially the same as those in
\cite{gutierrez-huang98}.

First we define a \emph{bowl-shaped domain} in spacetime and its
\emph{parabolic boundary.}  A set $\Omega\subset\R^n\times \R$ is
bowl-shaped if there are constants $t_0 < T$ so that $$\Omega =
\bigcup_{t_0\le t\le T} \Omega_t\times \{t\},$$ where each
$\Omega_t$ is convex and $\Omega_{t_1}\subset \Omega_{t_2}$ whenever
$t_1<t_2$.  The parabolic boundary of $\Omega$ is then
$\partial\Omega \setminus (\Omega_T\times\{T\}).$

\begin{pro} \label{gh-est}
Let $u$ be a smooth solution to (\ref{gen-eq}) which is convex in
$y$, and let $\Omega$ be a bowl-shaped domain in space-time
$\R^n\times\R$ so that $u=0$ on the parabolic boundary of
$\Omega$. Let $\beta$ be any unit direction in space.

Then at the maximum point $P$ of the function $$w = |u|\,\,
\partial^2_{\beta\beta} u\,\, e^{\frac12(\partial_\beta u)^2},$$
$w$ is bounded  by a constant depending on only $\alpha$, $\rho$,
$u(P)$, $\nabla u(P)$ and $n$.
\end{pro}

\begin{proof}
Choose coordinates so that $\beta = (1,0,\dots,0)$ and so that at
a maximum point $P$ of $w$, $u_{ij} = \frac{\partial^2 u}{\partial
y^i \partial y^j}$ is diagonal (in order to bound all second
derivatives $u_{\beta\beta}$, it suffices to focus only on the
eigendirections of the Hessian of $u$).

Since $w$ is positive in $\Omega$ and $0$ on the parabolic
boundary, there is a point $P$ outside the parabolic boundary of
$\Omega$ at which $w$ assumes its maximum value.  We work with
$\log w$ instead of $w$.  Then at $P$,
$$(\log w)_i = 0, \qquad (\log w)_t \ge 0,
\qquad (\log w)_{ij} \le 0.$$
 Here we use $i,j,t$ subscripts for
partial derivatives in $y^i$, $y^j$ and $t$, and the last inequality
is as a symmetric matrix.  These equations become, at $P$,
\begin{eqnarray}
&\displaystyle\frac{u_i}u + \frac{u_{11i}}{u_{11}} + u_1 u_{1i} = 0,&
\label{grad-w-zero}\\
&\displaystyle\frac{u_t}u + \frac{u_{11t}}{u_{11}} + u_1 u_{1t} \ge 0,&
\label{t-deriv}\\
&\displaystyle \frac{u_{ij}}u - \frac{u_i u_j}{u^2} +
\frac{u_{11ij}}{u_{11}} - \frac{u_{11i}u_{11j}}{u_{11}^2} +
u_{1j}u_{1i} + u_1u_{1ij} \le 0.& \label{hess-neg}
\end{eqnarray}
To use (\ref{t-deriv}), we compute
\begin{eqnarray*}
u_{1t} &=& \left[ -\rho (\det u_{ij})^{-\alpha} \right]_1 \\
&=& (\det u_{ij})^{-\alpha} ( -\rho_1 + \alpha\rho\,
u^{ij}u_{ij1}),\\
 u_{11t}
&=& (\det u_{ij})^{-\alpha} \\
&&{}\times\left[ 2\alpha\rho_1 \,u^{ij}u_{ij1} -
\alpha^2\rho(u^{ij}u_{ij1})^2 - \rho_{11} - \alpha \rho\,
u^{ik}u^{jl}u_{kl1}u_{ij1} + \alpha \rho\, u^{ij}u_{ij11}\right],
\end{eqnarray*}
where $u^{ij}$ is the inverse matrix of the Hessian $u_{ij}$. Now
plug into (\ref{t-deriv}) and divide out by $(\det
u_{ij})^{-\alpha}$ to find
 \begin{eqnarray} \nonumber
 &\displaystyle
\frac1{u_{11}} \left[2\alpha\rho_1 \,u^{ij}u_{ij1} -
\alpha^2\rho(u^{ij}u_{ij1})^2 - \rho_{11} - \alpha \rho\,
u^{ik}u^{jl}u_{kl1}u_{ij1} + \alpha \rho\, u^{ij}u_{ij11} \right]&
\\ &\displaystyle {}-\frac\rho u + u_1 ( -\rho_1 + \alpha\rho\,
u^{ij}u_{ij1}) \ge 0 & \label{u-11t}
\end{eqnarray}

The last term of the first line of (\ref{u-11t}) leads us to
contract (\ref{hess-neg}) with the positive-definite matrix $u^{ij}$
so that at $P$:
\begin{eqnarray*}
0 &\ge& u^{ij}\left(\frac{u_{ij}}u - \frac{u_i u_j}{u^2} +
\frac{u_{11ij}}{u_{11}} - \frac{u_{11i}u_{11j}}{u_{11}^2} +
u_{1j}u_{1i} + u_1u_{1ij} \right) \\
&=& \frac nu -\frac{2u^{ij}u_iu_j}{u^2} +
\frac{u^{ij}u_{11ij}}{u_{11}} - \frac{u^{ij}u_iu_1u_{1j}}u -
\frac{u^{ij}u_ju_1u_{1i}}u \\
&&{}- u^{ij}u_1^2u_{1i}u_{1j} + u^{ij}u_{1j}u_{1i} +
u^{ij}u_1u_{1ij} \qquad \qquad \mbox{(by (\ref{grad-w-zero}))}\\
&=& \frac nu -\frac{2u^{ij}u_iu_j}{u^2} +
\frac{u^{ij}u_{11ij}}{u_{11}} - \frac{2u_1^2}u  - u_1^2u_{11} +
u_{11} + u^{ij}u_1u_{1ij} \\
&& \qquad \qquad \mbox{(since $u_{ij}$ is diagonal at $P$)}\\
 &\ge& \frac nu -\frac{2u^{ij}u_iu_j}{u^2} - \frac{2u_1^2}u
- u_1^2u_{11} + u_{11} + u^{ij}u_1u_{1ij} + \frac1{\alpha u} +
\frac{\rho_1 u_1}{\alpha \rho} \\
&&{}- u_1u^{ij}u_{ij1} - \frac{2\rho_1 u^{ij}u_{ij1}}{\rho u_{11}} +
\frac{\alpha (u^{ij}u_{ij1})^2}{u_{11}} + \frac{\rho_{11}}{\alpha
\rho u_{11}} + \frac{u^{ik}u^{jl}u_{kl1}u_{ij1}}{u_{11}} \\
&& \qquad \qquad \mbox{(by (\ref{u-11t}))}\\
 &\ge & \frac
{n+\frac1\alpha}u - 2 \sum_{i=1}^n\frac{u_i^2}{u^2u_{ii}} -
\frac{2u_1^2}u - u_1^2u_{11} + u_{11} +
\frac{\rho_1 u_1}{\alpha \rho} \\
&&{} - \frac{\rho_1^2}{\alpha\rho^2 u_{11}} +
\frac{\rho_{11}}{\alpha \rho u_{11}} + \sum_{i,j=1}^n
\frac{u_{ij1}^2}{u_{11}u_{ii}u_{jj}}
\end{eqnarray*}
by collecting terms, completing the square, and since $u_{ij}$ is
diagonal at $P$. Continue computing
\begin{eqnarray*}
0 &\ge& \frac {n+\frac1\alpha}u -
2\sum_{i=1}^n\frac{u_i^2}{u^2u_{ii}} - \frac{2u_1^2}u - u_1^2u_{11}
+ u_{11} +
\frac{\rho_1 u_1}{\alpha \rho} \\
&&{} - \frac{\rho_1^2}{\alpha\rho^2 u_{11}} +
\frac{\rho_{11}}{\alpha \rho u_{11}}  + \frac{u_{111}^2}{u_{11}^3} +
2\sum_{i=2}^n \frac{u_{11i}^2}{u_{11}^2u_{ii}} \\
&=& \frac {n+\frac1\alpha}u - \frac{2u_1^2}{u^2u_{11}} -
\frac{2u_1^2}u - u_1^2u_{11} + u_{11} +
\frac{\rho_1 u_1}{\alpha \rho} \\
&&{} - \frac{\rho_1^2}{\alpha\rho^2 u_{11}} +
\frac{\rho_{11}}{\alpha \rho u_{11}}  + \frac{u_1^2}{u_{11}u^2} +
\frac{2u_1^2}u + u_1^2u_{11}
\end{eqnarray*}
by (\ref{grad-w-zero}) and since $u_{ij}$ is diagonal at $P$.
Finally, collect terms so that
 $$0\ge
u_{11} + \left(\frac{n+\frac1\alpha}u + \frac{\rho_1u_1}{\alpha\rho}
\right) + \frac1{u_{11}} \left(-\frac{u_1^2}{u^2} -
\frac{\rho_1^2}{\alpha\rho^2} + \frac{\rho_{11}}{\alpha\rho}\right)
$$
and multiply each side of the inequality by $u^2u_{11}e^{u_1^2}$ to
find a quadratic inequality $$w^2 + a w + b \le 0$$ for
$w=|u|u_{11}e^{\frac12 u_1^2}$ at $P$ the point in $\Omega$ at which
the maximum of $w$ is achieved. The coefficients $a$ and $b$ involve
only $n$, $\alpha$, $\rho$, $u(P)$ and $u_1(P)$, and so there is an
upper bound of $w$ on $\Omega$ depending on only these quantities.
\end{proof}

This bounds $\partial^2_{ij}s$ away from infinity, which, together
with Andrews's speed estimate, shows that the ellipticity is
locally uniformly controlled in the interior of appropriate
bowl-shaped domains. In the next section, we use barriers
constructed from Calabi's example to find appropriate bowl-shaped
domains.

\section{Applying Guti\'errez-Huang's Estimate}
\label{upper-interior-barrier}

In this section, we find bowl-shaped domains which are uniformly large
on any compact subset of $\mathcal D^\circ(s)$ for $s=s(\cdot,t)$ the
support function of $\mathcal K(t)$.  Upper barriers will be produced
from Calabi's example to achieve this.

First, we give an outline of our approach: Given a sequence of
smooth solutions $s$ to $\partial_ts = -(\det
s_{ij})^{-\frac1{n+2}}$, and a point $y_0$ in $\mathcal D^\circ
(s_0)$, modify $s$ by adding a linear function so that $s|_{t=0}$
has its minimum at $y_0$.  Adding a linear function does not affect
the flow. Then, if $s(y_0,0)=p$, the sublevel set $\{(y,t) :
s(y,t)\le p, 0\le t\le T\}$ is a bowl-shaped domain with $(y_0,0)$
at its vertex.  In order to apply Guti\'errez-Huang's estimate, we
must ensure that these bowl-shaped domains are uniformly large. This
amounts to showing that $s$ must decrease by a definite amount in a
neighborhood of $(y_0,0)$.

We achieve this by using barriers made out of Calabi's Example
\ref{calabi-ex}.  For each of the $n+1$ faces in Calabi's initial
orthant $\mathcal C=\mathcal C(0)$, consider the outward normal
directions $Y^i$, $i=1,\dots,n+1$. Under the affine normal flow of
Calabi's example, the support function $s_{\mathcal C}(Y^i,t)$
remains constant in $t$ for each $i=1,\dots,n+1$.  Thus Calabi's
example, in and of itself, is inadequate as a barrier to move the
support function in directions normal to these faces.

For our initial convex domain $\mathcal K$, we will obtain estimates
only for those $Y\in \mathcal D^\circ(s_{\mathcal K})$ the interior
of the domain of the support function. For such a $Y$, there is a
supporting hyperplane to $\mathcal K$ with $Y$ as its outer normal
which intersects the hypersurface $\partial \mathcal K$ in a compact
set $W$.  Then an appropriate barrier can be constructed as the
intersection
$$X=\bigcap_{i=1}^{n+1} \mathcal C^i$$
of $n+1$ affine images $\mathcal C^i$ of Calabi's example so that
the boundary of $X$ has one bounded face which contains $W$ and is
normal to $Y$, and $n+1$ unbounded faces (for example, a U-shaped
well in $\R^2$ is the intersection of two affine images of the first
quadrant).  Under the affine normal flow, $X(t)$ must remain inside
each $\mathcal C^i(t)$, and the explicit solution to Calabi's
example then shows that the support function $s_{X}(Y)$ must move as
$t$ increases away from 0.  The discussion below proves this
geometric sketch by working with the support functions instead of
the hypersurfaces involved.

Here are the details of the construction. Recall
Calabi's Example \ref{calabi-ex} from above:
$$\mathcal C(t) = \{(x^1,\dots,x^{n+1})\in\R^{n+1} :
 x^i\ge0, \prod_i x^i  = c_n t^{\frac{n+2}2}\}$$
 for $c_n = (n+1)^\frac12 (\frac2{n+2})^{\frac{n+2}2}$ and $t\ge0$.
Compute the support function for $t\ge0$ and
$Y=(y^1,\dots,y^{n+1})$:
\begin{equation} \label{calabi-support}
s_{\mathcal C}(Y,t)=\left\{
\begin{array}{cl} +\infty & \mbox{ if any }y^i>0 \\
-(n+1)\left(c_nt^{\frac{n+2}n}\prod_{i=1}^{n+1}|y^i|
\right)^{\frac1{n+1}} & \mbox{ if all } y^i\le 0
\end{array} \right.
\end{equation}
In our analysis, we restrict the homogeneity-one function
$s_{\mathcal C}$ to an affine hyperplane $\sim\R^n$ so that
$s_{\mathcal C}(0)=0$ on a simplex in $\R^n$ and is $+\infty$
elsewhere.

Now consider the action of the affine group on the support function.
If $\mathcal K$ is a convex body and $Y\in\R^{n+1}$, recall $s_{\mathcal K}(Y) =
\sup_{x\in\mathcal K} \langle x,Y\rangle$. Then for any matrix $A$
and vector $b$,
 \begin{equation} \label{aff-ch-support}
 s_{A \mathcal K + b}(Y) = s_{\mathcal K}(Y)(A^\top
Y) + \langle b,Y \rangle.
 \end{equation}
Therefore, for any simplex $\mathcal S$ in $\R^n$ and any linear
function $\ell(y)$ on $\R^n$, there is an affine image of $\mathcal
C(0)$ whose support function restricted to an affine
$\R^n\subset\R^{n+1}$ is equal to $\ell(y)$ on its domain $\mathcal
S$.

We will use $n+1$ of these copies of Calabi's example to construct
our barrier.  Consider a regular $(n+1)$-simplex in $\R^{n+1}$ with one
vertex at the origin and so that the face opposite this vertex is in
a hyperplane $y^{n+1}=c>0$ and intersects the positive $y^{n+1}$
axis.  Then the $n+1$ remaining faces of this simplex form the graph
of a piecewise-linear convex function $P$ whose domain is a simplex
$\mathcal S_n$ in $\R^n$ centered at the origin.  Extend this
function to be $+\infty$ outside $\mathcal S_n$.  We refer to $P$ as
a \emph{polyhedral pencil} function, after the shape of the region
above its graph.

Now consider our convex body $\mathcal K = \cup_{m=1}^\infty
\mathcal K^i$, and let $s_{\mathcal K}(y)$ denote the support
function of $\mathcal K$ restricted to an affine slice of
$\R^{n+1}$.  Let $\mathcal N$ be a compact subset the interior of
the domain of $s_{\mathcal K}(y)$ (i.e.\ $\mathcal N$ is the
intersection of the affine hyperplane $\R^n$ with a compact subset
of $\mathcal D^\circ(s_{\mathcal K})$). Then we know that
$s_{\mathcal K^m} = s_m \to s = s_{\mathcal K}$ uniformly on
$\mathcal N$ and that $|ds_m|$ is uniformly bounded on $\mathcal
N$.  This bound on the first derivatives of $s_m$ means that there
is a uniform $\lambda>1$ so that by replacing $P(y)$ by
$\lambda^nP(\lambda y)$, we have
 \begin{equation} \label{pencil-barrier}
 P(y-\tilde y)+ \sum_{j=1}^n
\frac{\partial s_m}{\partial y^j}(\tilde y) (y^j-\tilde y^j) +
s_m(\tilde y)\ge s_m(y)
\end{equation}
for all $\tilde y\in \mathcal N$, $y\in \R^n$.

So the polyhedral pencil function $P$ provides an initial barrier
at each point $\tilde y\in\mathcal N$.  We do not have an explicit
solution for the evolution of $P$, but we can conclude enough to
apply Guti\'errez-Huang's estimates.  Since $P$ can be extended to
be a convex, lower-semicontinuous function of homogeneity one on
$\R^{n+1}$, there is a corresponding convex body $\mathcal K_P$
whose support function is $P$.  The affine normal flow on
$\mathcal K_P$ induces a natural flow on $P$: $P(Y,0)=P(Y)$ for
all $Y\in \R^{n+1}$, and $P(Y,t)$ is the support function of
$\mathcal K_P(t)$. Assume that the affine hyperplane $\R^n\subset
\R^{n+1}$ is given by $\{y^{n+1}=-1\}$, and then we have

\begin{lem} \label{P-decreases}
$P(0,\dots,0,-1,t)<0$ for all $t>0$.
\end{lem}

Note that in the notation $Y=(y,-1)$, $P(0,\dots,0,-1,t)$ is just
$P(y,t)$ for $y=0$. (So the Lemma may be restated as $P(0,t)<0$ for
all $t>0$.) We use this notation for the proof.

\begin{proof}
Note $P(0,0)=0$.

As a function of $y$, the domain of $P(y,0)$ consists of $n+1$
simplices in $\R^n$, and $P(y,0)$ is the restriction of a linear
function on each one.  $P(y,0)$ is then the minimum of $n+1$ copies
$\mathcal C_1,\dots \mathcal C_{n+1}$ of Calabi's initial example,
each properly modified by an affine transformation.  Each of
$\mathcal C_1(t), \dots \mathcal C_{n+1}(t)$ is an upper barrier for
the evolution of $P$.  $\mathcal C_k(0,t)=0$ for all $t\ge 0$,
however.

To show that $P(0,t)<0$, we use the fact that $P(y,t)$ is always
convex and less than each $\mathcal C_k(y,t)$.  The explicit formula
(\ref{calabi-support}), together with (\ref{aff-ch-support}), shows
that each $\mathcal C_k(y,t)<0$ for $y$ near zero on a ray $R_k$
leaving the origin---this is because, near the origin, the $-(\prod
|y^i|)^\frac1{n+1}$ term in (\ref{calabi-support}) will dominate any
linear term coming from (\ref{aff-ch-support}).  Since $P(y,t)$ is
convex in $y$ and is less than each $\mathcal C_k(y,t)$, the graph
of $P(y,t)$ must be below the convex hull of the graphs of
$\{\mathcal C_k(y,t) \}_{k=1}^{n+1}$.  Since 0 is in the convex hull
of the rays $\{R_k\}_{k=1}^{n+1}$ (because $P$ was constructed using
a regular $(n+1)$-simplex in $\R^{n+1}$) and since $\mathcal C_k(y,t)<0$ for
$y\in R_k$ near 0, we conclude $P(0,t)<0$ for each $t>0$.
\end{proof}

For each $\tilde y\in \mathcal N$, and for $m=1,2,3,\dots$,
consider
$$\tilde s_m(y) = s_m(y) - s_m(\tilde y) - \sum_{j=1}^n
\frac{\partial s_m}{\partial y^j}(\tilde y) (y^j-\tilde y^j).$$
Then at $t=0$, $\tilde s_m(y,0)$ has its minimum value of $0$ at
$y=\tilde y$. As time goes forward, for each $T>0$, the sublevel
set $\{(y,t): t\in(0,T],\, \tilde s_m(y,t)<0\}$ is a bowl-shaped
domain. This bowl-shaped domain must contain the sublevel set
$$\mathcal B = \{(y,t): t\in(0,T],\, P(y-\tilde y)<0\},$$ which
contains $\{0\}\times (0,T]$ by Lemma \ref{P-decreases}. Note that
$\mathcal B$ is (except for translation) independent of $m$ and
$\tilde y\in \mathcal N$.  There is an increasing, positive
function  of $t>0$ $\epsilon(t)$ so that for each $\tilde y \in
\mathcal N$, Guti\'errez-Huang's estimates can be applied
uniformly on the ball $B_{\epsilon(t)/2}(\tilde y)\times \{t\}$ to
$\tilde s_m$---since $\tilde s_m$ satisfies the same flow equation
(\ref{s-y-flow}) as $s_m$.

Since the second derivatives of $\tilde s_m$ are the same as those
of $s_m$, Guti\'errez-Huang's estimate Proposition \ref{gh-est},
Andrews's speed bound Proposition \ref{andrews-speed} and the
convexity of $s_m$ imply uniform $C^2$ estimates on $s_m$ on each
compact subset of $\mathcal D^\circ(s) \times (0,T]$, where $T$ is
chosen so that each $s_m\ge r$ on $[0,T]$ (this is possible for
some $T$ by choosing coordinates so that an evolving sphere
centered at the origin as a uniform inner barrier.)

\begin{pro} \label{loc-unif-C2}
If $T$ is chosen so that each $s_m\ge r$ on $\sph^n\times [0,T]$,
then on each compact subset of $\mathcal D^\circ(s) \times (0,T]$
there are uniform spatial $C^2$ estimates for $s_m$ and the Hessian
of $s_m$ is uniformly bounded away from zero.
\end{pro}

Recall that $s_{\mathcal K^m}\to s_{\mathcal K}$ everywhere on
$\R^{n+1}\times [0,T]$ by Propositions \ref{exhaust-flow} and
\ref{support-exhaust-limit}.  The estimates will give greater
regularity to this pointwise convergence.

Note that the locally uniform spatial $C^2$ estimates in Proposition
\ref{loc-unif-C2} imply, by the evolution equation (\ref{s-y-flow}),
locally uniform parabolic $C^{2,1}$ estimates (i.e.\ two derivatives
in spatial coordinates and 1 in $t$). Then, since the logarithm of
the Monge-Amp\`ere operator is concave, Krylov's interior parabolic
$C^{2+\alpha,1+\frac\alpha2}$ estimates \cite{krylov87} are
available. Ascoli-Arzel\'a then shows that the convergence must be
in $C^{2,1}$ on each compact subset of $\mathcal D^\circ(s)\times
(0,T]$. Indeed, $s_{\mathcal K}$ is a $C^{2+\alpha,1+\frac\alpha2}$
solution on $\mathcal D^\circ(s)\times (0,T]$, and further
bootstrapping shows $s_{\mathcal K}$ is smooth. See e.g.\
Guti\'errez-Huang \cite{gutierrez-huang98} for details on defining
the $C^{2+\alpha,1+\frac\alpha2}$ norm and on applying Krylov's
estimates in the present case.

A remaining issue is long-time regularity. Since long-time existence
is already guaranteed, we need only apply the estimates again
starting at $t=T$.  The only possible sticking point is that we
still need to make sure that the same $r$ satisfying $s_{\mathcal
K^m}\ge r$ still works (in order to apply Andrews's speed estimate
Proposition \ref{andrews-speed}). This can be assured by an affine
change of coordinates. As in the proof of Theorem
\ref{long-time-ex}, $\mathcal K$ contains ellipsoids of arbitrarily
large volume.  We can change the affine coordinates so that an
appropriate ellipsoid becomes a sphere centered at the origin which
is large enough to guarantee that $s_{\mathcal K}(Y,t)> 2r $ for all
$Y\in\sph^n$ and $t\in[T,2T]$. This is certainly enough to ensure
that we can choose new exhausting domains $\mathcal K^m$ satisfy
$s_{\mathcal K^m}(Y,t)\ge r$ for all $Y\in \sph^n$ and $t\in[T,2T]$.

\begin{thm} \label{support-smooth}
If $\mathcal K$ is an unbounded convex domain in $\R^{n+1}$ which
contains no lines, then, under the affine normal flow, the support
function $s_{\mathcal K}= s_{\mathcal K}(Y,t)$ is smooth and
spatially locally strictly convex on $\mathcal D^\circ(s_{\mathcal
K}) \times (0,\infty)$.
\end{thm}

\section{Regularity of the hypersurface} \label{hyp-reg-sec}

We've seen in the previous sections that under the affine normal
flow, if $\mathcal K$ is an unbounded convex domain in $\R^{n+1}$
containing no lines, the support function $s_{\mathcal K}$ evolves
to be smooth and strictly convex for all positive time for all
$Y\in \mathcal D^\circ(s_{\mathcal K})$.  In this section, we
verify that, for $t>0$, every supporting hyperplane of the
evolving hypersurface $\partial \mathcal K(t)$ has its normal
vector in $\mathcal D^\circ(s_{\mathcal K})$. Therefore, since the
smoothness and convexity of the hypersurface are determined by the
regularity of the support function, the hypersurface $\partial
\mathcal K(t)$ is smooth and strictly convex for all $t>0$.

\begin{thm} \label{smooth-hypersurface}
Let $\mathcal K$ be an unbounded convex domain in $\R^{n+1}$ which
contains no lines. Then, under the affine normal flow, the
hypersurface $\partial \mathcal K(t)$ is smooth and strictly convex
for all $t>0$.

Moreover, if
$$\mathcal K = \bigcup_{i} \mathcal
K^i, \qquad \mathcal K^{i} \subset \mathcal K^{i+1},$$ where each
$\mathcal K^i$ is bounded and have smooth, strictly convex boundary,
then for all $t>0$,  each $p\in\partial \mathcal K(t)$ has a
neighborhood on which the sequence of  hypersurfaces $\mathcal
\partial K^i(t)$ converges to $\partial \mathcal K(t)$ in the
$C^\infty$ topology.
\end{thm}

The proof will depend on finding appropriate initial barriers.  We
begin with some easy results on the support function.

\begin{lem} \label{bound-domain-support}
If $\mathcal K$ is an unbounded convex domain in $\R^{n+1}$, then
for every nonzero $Y_0\in\partial \mathcal D(s_{\mathcal K})$, there
is a ray $R$ perpendicular to $Y_0$ which is contained in the
closure $\mathcal{\bar K}$.
\end{lem}

\begin{proof}
We work in terms of support functions.  The support function
$s_{\mathcal K}$ is a homogeneity-one, convex,
lower-semicontinuous function on $\R^{n+1}$ with values in
$(-\infty,+\infty]$.  Since $\mathcal K$ is unbounded,
$s_{\mathcal K}$ must assume the value $+\infty$, and the
convexity of $s_{\mathcal K}$ implies $s_{\mathcal K}$ is infinite
on an open half-space of $\R^{n+1}$.

$R\subset \mathcal {\bar K}$ if and only if $s_R\le s_{\mathcal
K}$ on all of $\R^{n+1}$.  For the ray $$R = \{w+\tau v : \tau\ge
0\},\qquad s_R(Y) = \left\{ \begin{array}{c@{\mbox{ for }}c}
\langle Y,w\rangle & \langle Y,v \rangle \le 0 \\
+\infty & \langle Y,v\rangle > 0.
\end{array} \right.$$
Thus, given $s_{\mathcal K}$ and $Y_0\in \partial \mathcal D
(s_{\mathcal K})$, we seek an $R$ so that $R\perp Y_0$ and $s_R\le
s_{\mathcal K}$.

Since $\mathcal D(s_{\mathcal K})$ is a convex cone in $\R^{n+1}$
with vertex at the origin, if $Y_0 \in\partial \mathcal
D(s_{\mathcal K})$, then $\mathcal D(s_{\mathcal K})$ is contained
in a closed half-space with $Y_0$ in its boundary. Thus there is a
nonzero vector $v$ so that $$\mathcal D(s_{\mathcal K}) \subset \{Y:
\langle Y,v\rangle\le 0\}, \qquad \langle Y_0,v\rangle = 0.$$ In
order to find $R$, we also need a vector $w$ so that $\langle Y,w
\rangle \le s_{\mathcal K}(Y)$ for all $Y\in \mathcal D(s_{\mathcal
K})$.  This is easy: $\langle Y,w\rangle$ is the support function of
the convex set $\{w\}$. So for any $w\in \mathcal {\bar K}$,
$\langle Y,w\rangle \le s_{\mathcal K}(Y)$ for all $Y\in \R^{n+1}$,
and $R=\{v+\tau w : \tau\ge 0\}$ is the ray to be constructed.
\end{proof}

\begin{lem} \label{exist-tube}
If $\mathcal K$ is an unbounded convex domain in $\R^{n+1}$, $Y\in
\partial \mathcal D(s_{\mathcal K})$, and $R$ is any ray
contained in $\mathcal {\bar K}$, there is a half-cylinder $\mathcal
Q$ pointing in the direction of $R$ which is contained in the open
set $\mathcal K$.
\end{lem}

\begin{proof}
Let $B$ be an open ball contained in $\mathcal K$. Then the convex
hull of $R$ and $B$ contains such a half-cylinder.
\end{proof}

We are now ready to prove Theorem \ref{smooth-hypersurface}.

\begin{proof}[Proof of Theorem \ref{smooth-hypersurface}]
We show that for all $t>0$, every supporting hyperplane of $\mathcal
K(t)$  must have outward normal vector $Y_0$ lying in $\mathcal
D^\circ(s_{\mathcal K})$.  Then the smoothness and strict convexity
of the support function $s_{\mathcal K(t)}$  imply that the
hypersurface $\partial \mathcal K(t)$ is also smooth and strictly
convex.

First we rule out the case $Y_0\notin \overline{\mathcal
D(s_{\mathcal K})}$.  In this case, there is a closed half-space
of $\R^{n+1}$ containing $\mathcal D(s_{\mathcal K})$ but
excluding $Y_0$.  In other words, there is a nonzero vector $v$ so
that
$$\mathcal D(s_{\mathcal K}) \subset \{Y: \langle Y,v\rangle\le 0\},
\qquad \langle Y_0,v\rangle > 0.$$ Then, as in Lemmas
\ref{bound-domain-support} and \ref{exist-tube} above, there is a
half-cylinder $\mathcal Q$ in the direction of $v$ contained in
$\mathcal K$.  Since there are ellipsoids of arbitrarily large
volume inside $\mathcal Q$ to act as barriers, $\mathcal Q$ always
intersects $\mathcal K(t)$. Since $\mathcal Q$ is in the direction
of $v$ and $\langle Y_0,v\rangle >0$, this shows that $s_{\mathcal
K(t)}(Y_0)=+\infty$ for all $t>0$.  Since $\mathcal K(t)$ is convex,
this shows it has no supporting hyperplane with outward normal
$Y_0\notin \overline{\mathcal D (s_{\mathcal K})}$.

Finally, we show that if $Y_0\in \partial \mathcal D(s_{\mathcal
K})$ is a nonzero vector, then there is no supporting hyperplane to
$\partial \mathcal K(t)$ with outward normal $Y_0$. By Proposition
\ref{boundary-fixed} below, $s_{\mathcal K(t)}(Y_0)= s_{\mathcal
K}(Y_0)$ for all $t>0$.  Thus, we simply need to ensure that the
hyperplane $\mathcal P = \{x: \langle Y_0,x \rangle = s_{\mathcal
K}(Y_0)\}$ does not intersect $\overline{\mathcal K(t)}$ for $t>0$.
To achieve this, we choose an affine image $\mathcal I$ of Calabi's
example as an initial outer barrier. In particular, one of the faces
of $\mathcal I$ can be chosen to be contained in the hyperplane
$\mathcal P$. (Proof: The support function of $\mathcal I$ is a
linear function on a cone over a simplex and is $+\infty$ elsewhere.
To find such a function to be an upper barrier to $s_{\mathcal K}$
at $Y_0$, note that for any closed simplex contained in $\mathcal
D(s_{\mathcal K})$, the support function $s_{\mathcal K}$ is
continuous on this simplex by Theorem 10.2 in \cite{rockafellar}. So
for any cone $\mathcal C$ over a closed $n$-simplex containing $Y_0$
and contained in $\mathcal D(s_{\mathcal K})$, we may find a linear
function as an upper barrier to $s_{\mathcal K}$ at $Y_0$. Then
extend this function to be $+\infty$ outside $\mathcal C$.) The
explicit solution to Calabi's example proves that $\mathcal P$ does
not intersect $\overline{\mathcal I(t)} \supset \overline{\mathcal
K(t)}$ for all $t>0$.

Thus all the supporting hyperplanes of $\partial \mathcal K(t)$ have
outward normal in $\mathcal D^\circ (s_{\mathcal K})$, and Theorem
\ref{smooth-hypersurface} is proved.
\end{proof}

\begin{pro} \label{boundary-fixed}
If $\mathcal K$ is a convex unbounded domain in $\R^{n+1}$ which
does not contain any lines, then under the affine normal flow, the
support function $s_{\mathcal K(t)}(Y_0) = s_{\mathcal K}(Y_0)$ for
all $t>0$ and $Y_0\in \partial \mathcal D(s_{\mathcal K})$.
\end{pro}

\begin{proof}
It is obvious that $s_{\mathcal K(t)}(Y_0) \le s_{\mathcal K}
(Y_0)$ since the effect of the affine normal flow on support
functions is to decrease them.  We need only show $s_{\mathcal
K(t)}(Y_0) \ge s_{\mathcal K} (Y_0)$ for $Y_0\in \partial \mathcal
D(s_{\mathcal K})$.

We achieve this by using ellipsoids as inner barriers.  Assume
that $s_{\mathcal K}(Y_0)<+\infty$ and let $\epsilon>0$.  Then
there is an $x\in\mathcal K$ so that $$\langle x,Y_0 \rangle >
s_{\mathcal K}(Y_0) - \epsilon.$$ Lemmas
\ref{bound-domain-support} and \ref{exist-tube} ensure that there
is a half-cylinder $\mathcal Q\subset \mathcal K$ which points in
a direction $v$ perpendicular to $Y_0$. Then, inside the convex
hull of $\mathcal Q$ and $\{x\}$, there is another half-cylinder
$\mathcal Q'\subset \mathcal K$ which points in the direction of
$v$ and whose central ray contains a point $x'$ satisfying
$$\langle x',Y_0 \rangle > s_{\mathcal K}(Y_0) - 2\epsilon.$$  Now
there are ellipsoids of arbitrarily large volume contained in
$\mathcal Q'$, and these inner barriers show that for all $t>0$,
there is a point $x''$ on the central ray of $\mathcal Q'$ which is
contained in $\mathcal K(t)$. Now since $x''-x'$ is perpendicular to
$Y_0$,
$$ s_{\mathcal K(t)}(Y_0) \ge \langle x'',Y_0 \rangle = \langle
x',Y_0 \rangle > s_{\mathcal K}(Y_0) - 2\epsilon.$$ Thus
$s_{\mathcal K(t)}(Y_0) \ge s_{\mathcal K}(Y_0)$ so long as
$s_{\mathcal K}(Y_0)<+\infty$. The case $s_{\mathcal
K}(Y_0)=+\infty$ is essentially the same.
\end{proof}

\section{A Dirichlet Problem} \label{dirichlet-section}

Proposition \ref{boundary-fixed} above shows that the affine normal
flow on noncompact domains can be recast as a Dirichlet boundary
problem for the support function, although discontinuous and
infinite boundary values are allowed. In the interior $\mathcal
D^\circ (s_{\mathcal K})$, the support function evolves by the
affine normal flow equation, while the value of the support function
on the boundary $\partial \mathcal D(s_{\mathcal K})$ is fixed.  At
each positive time $t$, $s_{\mathcal K(t)}$ is lower-semicontinuous.


In terms of PDEs, we can take an affine slice of the domain of the
support function.  Consider first the case when $\mathcal
D(s_{\mathcal K})$ contains no lines (this is true if and only if
$\mathcal K$ contains a nonempty open convex cone).  In this case,
we can choose coordinates so that $\Omega = \{y\in\R^n : (y,-1)\in
\mathcal D^\circ(s_{\mathcal K})\}$ is bounded. The support
function, when restricted to this hyperplane, then satisfies the
Dirichlet boundary problem for the flow equation
$$ \frac{\partial s}{\partial t} = - \left( \det \frac{\partial^2
s}{\partial y^i \partial y^j} \right)^{-\frac1{n+2}}$$ with initial
condition given by $s_{\mathcal K}$.   If $\mathcal
D^\circ(s_{\mathcal K})$ does contain a line, then we must consider
more than one affine hyperplane slice.  Since $s_{\mathcal K}$ has
homogeneity one, this amounts to considering $s_{\mathcal K}$ as a
section of the tautological bundle over projective sphere $\sph^n_P
= (\R^{n+1}\setminus \{0\})/ \R^+$, where $\R^+$ acts on $\R^{n+1}$
by homothetic scaling.

Alternately, we can consider $s=s_{\mathcal K}$ restricted to the
Euclidean sphere $\sph^n$.  Define a subset of $\sph^n$ to be convex
if it is the intersection of $\sph^n$ with a convex cone in
$\R^{n+1}$ with vertex at the origin.  Then $\Upsilon = \sph^n \cap
\mathcal D^\circ (s)$ is a convex domain in $\sph^n$, and $s$
evolves under the affine normal flow via a Dirichlet problem on
$\Upsilon$ with equation, as in \cite{andrews00}, for
$s=s|_{\sph^n}$
$$ s_t  = - \left[ \det \left( s_{;ab} + s\delta_{ab}\right)\right]
^{-\frac1{n+2}}.$$ Here $s_{;ab}$ denotes second covariant
derivative of $s$ with respect to the standard connection on
$\sph^n$ and the subscripts $a,b$ indicate an orthonormal frame.

It is an interesting question to study under what condition this
Dirichlet problem admits a unique solution. We plan to study this
problem in detail later. We remark now that in the case that when
$s$ is continuous when restricted to the boundary $\partial \mathcal
D(s)$, then $s$ must be continuous on the closure
$\overline{\mathcal D(s)}$ (see Lemma \ref{cont-bound} below).  Thus
in the case $s$ is continuous and finite when restricted to
$\partial D(s)$, the Dirichlet problem has a unique solution by the
maximum principle.

\begin{lem} \label{cont-bound} Let $s$ be a convex, lower
semicontinuous function from $\R^n$ to $(-\infty,\infty]$. If $s$ is
continuous when restricted to $\partial D(s)$, then it is continuous
on the closure of its domain $\overline{\mathcal D(s)}$.
\end{lem}

\begin{proof}
Let $x_i\in \mathcal D^\circ(s)$, $x_i\to x\in\partial \mathcal
D(s)$. Let $z\in\mathcal D^\circ(s)$ and let $y_i$ be the
intersection of $\partial D(s)$ and the ray from $z$ to $x_i$.
$y_i\to x$ and so $s(y_i)\to s(x)$. Moreover, $s$ is convex
restricted to each such ray, and so
$$s(y_i) - s(x_i) \ge \frac{|y_i-x_i|}{|x_i-z|} [s(x_i)-s(z)].$$ Thus,
since ${|y_i-x_i|}/{|x_i-z|}\to 0$,  $$\limsup s(x_i) \le \lim
s(y_i) = s(x).$$ Lower semicontinuity then shows $\lim s(x_i)=s(x)$.
\end{proof}

\section{Proofs of Theorems} \label{final-sec}

Here we restate Theorem \ref{classify-ancient} a bit more precisely:
\begin{thm}
Let $\mathcal L(t)$ be a solution to the affine normal flow defined
for all $t\in(-\infty,0)$.  Assume that at some $t_0\in(-\infty,0)$,
the convex hull $\widehat{\mathcal L(t_0)}$ contains no lines. Then
$\mathcal L(t)$ must be a paraboloid translating in time or an
ellipsoid shrinking in time.
\end{thm}

\begin{proof}
Consider $\mathcal L'(t)$ defined for $t\in[\tau,0)$ so that
$\mathcal L'(\tau)$ is smooth, compact, and strictly convex.  Then
Proposition \ref{c2-decay} and the semigroup property show that the
cubic form
$$|C|^2_{\mathcal L'(t)} \le \frac{c_n}{t-\tau}$$ for all
$t\in[\tau,0)$ for $c_n$ a constant depending only on the dimension.

Theorem \ref{smooth-hypersurface} then shows that $\mathcal L(t)$,
for $t\in[\tau,0)$, is locally a $C^\infty$ limit of such $\mathcal
L'(t)$. Thus $|C|^2_{\mathcal L(t)} \le c_n / (t-\tau)$ also.  Since
$\mathcal L(t)$ is an ancient solution we may let $\tau\to-\infty$.
So  $C_{ijk}=0$ identically on $\mathcal L(t)$ for all $t$.

A well-known classical theorem of Berwald (see e.g.\ Cheng-Yau
\cite{cheng-yau86} or Nomizu-Sasaki \cite{nomizu-sasaki}) shows that
$\mathcal L(t)$ must be a hyperquadric: an ellipsoid, a paraboloid,
or a hyperboloid.  Only the ellipsoid (a shrinking soliton) and the
paraboloid (a translating soliton) are part of an ancient solution.
\end{proof}






We can also prove the following existence result on hyperbolic
affine spheres which is essentially due to Cheng-Yau
\cite{cheng-yau77}. The essential step, due to Cheng-Yau, is to
solve the Monge-Amp\`ere equation
$$\det\partial^2_{ij}\phi = \left(-\frac1\phi\right)^{\frac1{n+2}}, \qquad
\phi|_{\partial\Omega} = 0, \qquad \partial^2_{ij}\phi>0$$ on a
convex bounded domain $\Omega\subset\R^n$. The radial graph of
$-\frac1\phi$ over $\Omega$ is then a hyperbolic affine sphere
asymptotic to the boundary of the cone over $\Omega$. We note that
the proper embeddedness of the hyperbolic affine sphere is contained
in Gigena \cite{gigena81} and Sasaki \cite{sasaki80}. See
\cite{loftin01} for a more complete discussion.

\begin{thm}
For every open convex cone $\mathcal C$ in $\R^{n+1}$ which contains
no lines, there is a properly embedded hyperbolic affine sphere in
$\R^{n+1}$ asymptotic to the boundary of $\mathcal C$.
\end{thm}

\begin{proof}
Example \ref{weak-hyp-as} ensures that under the affine normal flow,
the boundary $\partial\mathcal C$ evolves as an expanding soliton
$\partial \mathcal C(t)$. The regularity result Theorem
\ref{smooth-hypersurface} ensures that for $t>0$ the hypersurface
$\partial \mathcal C(t)$ is smooth and strictly convex. Thus, for
each $t>0$, $\partial \mathcal C(t)$ is a hyperbolic affine sphere
by Proposition \ref{soliton-affine-sphere}. The discussion in
Section \ref{dirichlet-section} shows that  $\partial \mathcal C(t)$
is asymptotic to the boundary of the cone $\mathcal C$.
\end{proof}

\bibliographystyle{abbrv}
\bibliography{thesis}

\end{document}